\documentclass[a4,10pt]{article}

\usepackage{algpseudocode}
\usepackage{algorithm,algcompatible} 

\algnewcommand\algorithmicinput{\textbf{Input:}}
\algnewcommand\INPUT{\item[\algorithmicinput]}
\algnewcommand\algorithmicoutput{\textbf{Output:}}
\algnewcommand\OUTPUT{\item[\algorithmicoutput]}
\usepackage{amsmath}
\usepackage{mathtools}
\usepackage{amssymb}
\usepackage{amsthm}
\usepackage[american]{babel}%
\usepackage{bm}
\usepackage{bbm}
\usepackage{booktabs}
\usepackage[top=30truemm,bottom=30truemm,left=25truemm,right=25truemm]{geometry}
\usepackage[bookmarksnumbered=true]{hyperref}
\hypersetup{
     colorlinks = true,
     linkcolor = blue,
     anchorcolor = blue,
     citecolor = blue,
     filecolor = blue,
     urlcolor = blue
}
\usepackage{comment}
\usepackage{color} 
\usepackage{csquotes}
\usepackage{enumitem}
\usepackage{epstopdf}
\usepackage{eso-pic}%
\usepackage{etoolbox}%
\usepackage{fancyvrb}
\usepackage[T1]{fontenc}%
\usepackage{graphicx}
\usepackage{grffile}%
\usepackage[
hyperref,%
table%
]{xcolor}%
\usepackage{ifthen}
\newboolean{PrintVersion}
\setboolean{PrintVersion}{false}
\ifthenelse{\boolean{PrintVersion}}{   
\hypersetup{  
    citecolor=black,%
    filecolor=black,%
    linkcolor=black,%
    urlcolor=black}
}{} 

\usepackage{authblk}
\usepackage{latexsym,indentfirst}
\usepackage{listings}
\usepackage{lmodern}%
\usepackage{lscape} 
\usepackage{multirow}%
\usepackage{microtype}%
\usepackage[round]{natbib}
\usepackage{rotating}%
\usepackage{siunitx}%
\usepackage{subfigure}
\usepackage{tabularx}%
\usepackage{threeparttable} 
\usepackage{tikz}
\usetikzlibrary{positioning}
\usetikzlibrary{calc}
\usetikzlibrary{shapes.geometric}
\usetikzlibrary{arrows.meta}
\usepackage{wrapfig}%

\theoremstyle{definition}
\newtheorem{theorem}{Theorem}[section]
\newtheorem{corollary}[theorem]{Corollary}

\newtheorem{proposition}[theorem]{Proposition}

\newtheorem{definition}[theorem]{Definition}
\newtheorem{example}[theorem]{Example}

\newtheorem{remark}[theorem]{Remark}

\xdefinecolor{lightgray}{RGB}{247, 247, 247}%
\xdefinecolor{semilightgray}{RGB}{240, 240, 240}%
\xdefinecolor{middlegray}{RGB}{127, 127, 127}%
\xdefinecolor{blue}{RGB}{58, 95, 205}%
\xdefinecolor{deepskyblue}{RGB}{0, 154, 205}%
\xdefinecolor{chocolate}{RGB}{205, 102, 29}%

\newcommand*{\N}{\operatorname{N}}

\newcommand*{\IN}{\mathbb{N}}

\newcommand*{\IR}{\mathbb{R}}

\newcommand{\argmax}{\mathop{\rm argmax}\limits}

\newcommand*{\E}{\mathbb{E}}

\newcommand*{\Prob}{\mathbb{P}}

\newcommand*{\rd}{\;\mathrm{d}}

\newcommand*{\Var}{\operatorname{Var}}

\newcommand{\bx}{\bm{x}}

\newcommand{\bone}{\bm{1}}

\newcommand{\darrow}{\stackrel{d}{\longrightarrow}}

\renewcommand*{\i}{{-1}}

\newcommand{\ou}[3]{%
  \mathrel{%
    \vcenter{\offinterlineskip
      \ialign{##\cr$#1$\cr\noalign{\kern-#3}$#2$\cr}%
    }%
  }%
}

\newcommand{\sqc}[2]{#1_{1},\dots,#1_{#2}}

\newcommand*{\T}{^{\top}}

\usepackage[normalem]{ulem}

\title{Measuring non-exchangeable tail dependence using tail copulas}

\author[1]{Takaaki Koike\footnote{Corresponding author; e-mail \url{takaaki.koike@r.hit-u.ac.jp}}}
\author[2]{Shogo Kato}
\author[3]{Marius Hofert}

\affil[1]{Graduate School of Economics, Hitotsubashi University, Kunitachi, Tokyo, 186-8601, Japan}
\affil[2]{Risk Analysis Research Center, Institute of Statistical Mathematics, Tachikawa, Tokyo, 190-8562, Japan}
\affil[3]{
Department of Statistics and Actuarial Science,
Faculty of Science,
The University of Hong Kong,
Pokfulam, Hong Kong}

\begin{document}
\maketitle
\date{}

\begin{abstract}
Quantifying tail dependence is an important issue in insurance and risk management.
The prevalent tail dependence coefficient (TDC), however, is known to underestimate the degree of tail dependence and it does not capture non-exchangeable tail dependence since it evaluates the limiting tail probability only along the main diagonal.
To overcome these issues, two novel tail dependence measures called the maximal tail concordance measure (MTCM) and the average tail concordance measure (ATCM) are proposed.
Both measures are constructed based on tail copulas and possess clear probabilistic interpretations in that the MTCM evaluates the largest limiting probability among all comparable rectangles in the tail, and the ATCM is a normalized average of these limiting probabilities.
In contrast to the TDC, the proposed measures can capture non-exchangeable tail dependence.
Analytical forms of the proposed measures are also derived for various copulas.
A real data analysis reveals striking tail dependence and tail non-exchangeability of the return series of stock indices, particularly in periods of financial distress.
\end{abstract}

\hspace{0mm}\\
\noindent \emph{Keywords:} 
Copula;
tail copula;
tail dependence; 
tail dependence coefficient;
tail dependence function;
tail non-exchangeability
\\
\noindent \emph{JEL codes:}
C01, 
C02, 
C40, 
G11. 

\section{Introduction}\label{sec:introduction}
The dependence between two continuous random variables $X\sim F$ and $Y\sim G$
is characterized by their \emph{copula} $C:[0,1]^2\rightarrow [0,1]$, that is
the distribution function of $(F(X),G(Y))$.  Of particular interest in extreme
value analysis is to quantify dependence in tail regions, so to summarize the
tendency for $X$ and $Y$ to jointly take on extremely small (or large) values.
According to \cite{ledford1996statistics,ledford1997modelling},
\cite{ramos2009new} and \cite{hua2011tail}, (lower) tail dependence of the
bivariate random vector $(X,Y)$ can be described by the \emph{tail order}
$1/\eta$, where $\eta\in (0,1]$, and the \emph{tail dependence parameter}
$\lambda\geq 0$ such that
\begin{align*}
C(p,p) \simeq l(p)\,p^{1/\eta}\quad (p \downarrow 0)\quad\text{and}\quad
\lambda=\lim_{p\downarrow 0}\,l(p)
\end{align*}
for some slowly varying function $l:\IR_{+}\rightarrow \IR_{+}=[0,\infty)$, where $f\simeq g$ ($x \rightarrow y$) for $f,g:\IR\rightarrow \IR$ means that $\lim_{x\rightarrow y}f(x)/g(x)=1$ for $y \in \IR\cup\{\pm \infty\}$.

The case when $\eta=1$ is of particular importance in insurance and risk management, for example, where the
\emph{Gaussian copula}
has been blamed
as a result of financial crisis of 2007--2009; see~\citet{embrechts2009linear} and \citet{donnelly2010devil}.  When
$\eta=1$, the tail dependence parameter
$\lambda(C)=\lim_{p\downarrow 0}C(p,p)/p$ is also known as the \emph{tail
  dependence coefficient}~\cite[TDC,][]{sibuya1960bivariate}.
  It is used for measuring tail dependence; see, for example,~\citet{aloui2011global} and~\citet{garcia2011dependence} for financial applications.
  Despite its
popularity, the TDC is known to underestimate the degree of tail dependence since
it quantifies the speed of decay of the joint tail probability only along the
main diagonal of $C$.
In addition, the TDC does not capture non-exchangeable tail dependence in the sense that $\lambda(C)$ always equals $\lambda(C\T)$, where $C\T$ is the copula of $(Y,X)$; see~\citet{hua2019assessing} and~\citet{bormann2020detecting} for recent analyses of tail non-exchangeability.
\cite{furman2015paths} addressed these issues and proposed
variants of the TDC where the main diagonal is replaced by the path maximizing
the joint tail probability.
Calculation and estimation of such tail indices may not always be straightforward due to the difficulty of deriving the path of maximal dependence for a given copula $C$; see \cite{sun2020statistical,sun2022tail} for recent progress on estimating such tail indices.
A similar measure of non-exchangeable tail dependence $\operatorname{limsup}_{(u,v)\downarrow(0,0)}C(u,v)/(u+v)$
 has also been considered in~\citet{genest2021class}.

\begin{figure}[t]
  \centering
  \includegraphics[width=15cm]{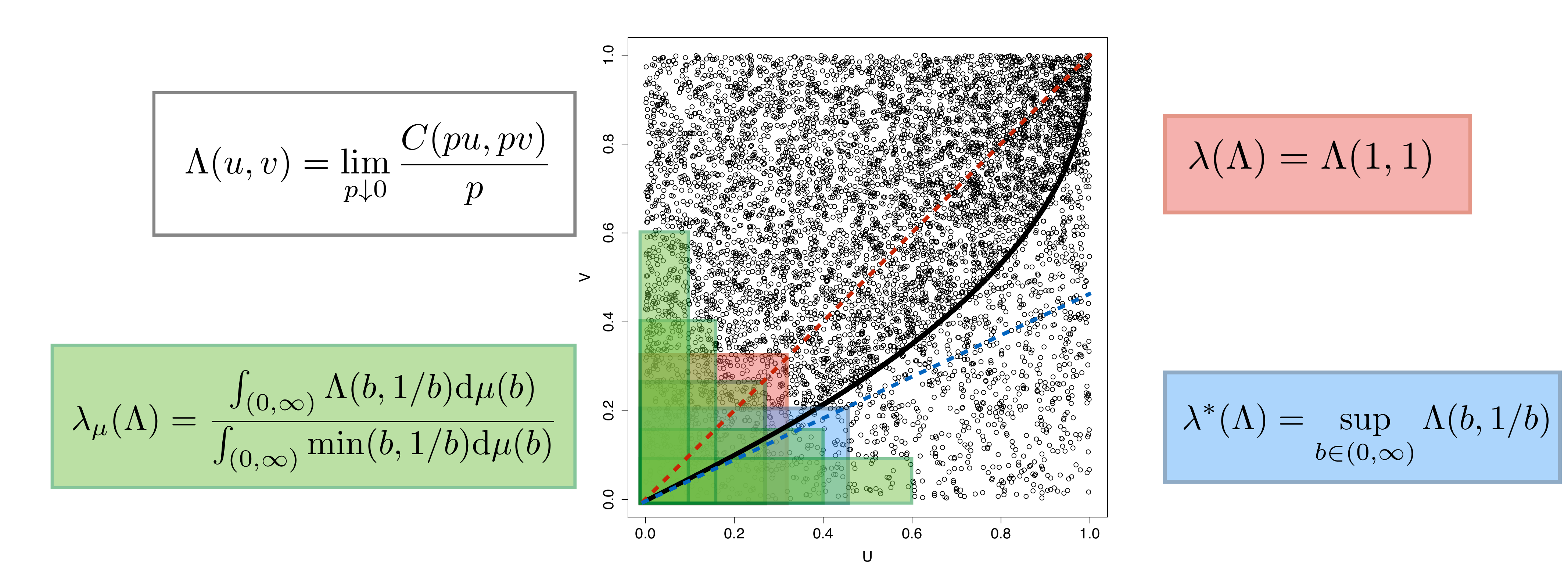}
  \caption{Scatter plot of the survival Marshall--Olkin copula
    $\hat C_{\alpha,\beta}^{\text{MO}}$ for $(\alpha,\beta)=(0.353,0.75)$, and
    its tail concordance measures.  The lower TDC $\lambda$ evaluates the limiting tail probability at the square colored in red.
    The $\mu$-ATCM $\lambda_\mu$ is a normalized average of $\Lambda$ over all comparable rectangles (green) weighted by the measure $\mu$.
    Finally, the MTCM $\lambda^\ast$ is the maximum of limiting tail probabilities over
    all comparable rectangles, which is attained for the rectangle highlighted in blue.
       }
  \label{fig:explanations:tail:dependence:measures}
\end{figure}

In this paper, we construct measures of non-exchangeable tail dependence, called the \emph{tail concordance measures (TCMs)}, based on the so-called \emph{tail copula}.
The tail copula of $(U,V)\sim C$ is defined by
\begin{align*}
\Lambda(u,v)
=\lim_{p\downarrow 0}\frac{C(pu,pv)}{p}=\lim_{p\downarrow 0}\frac{\Prob((U,V)/p \in [0,u]\times[0,v])}{p}
,\quad (u,v) \in \IR_{+}^2.
\end{align*}
The tail copula is a generalization of the TDC and plays an important role in extreme value analysis; see \cite{jaworski2004uniform}, \cite{schmidt2006non}, \cite{kluppelberg2007estimating}, \cite{nikoloulopoulos2009extreme} and \cite{joe2010tail}.
We construct measures of tail dependence based on the subset
\begin{align}\label{eq:ref:reference:set:Lambda:introduction}
\left\{\Lambda\left(b,\frac{1}{b}\right): b \in (0,\infty)\right\}.
\end{align}
Elements in the reference set~\eqref{eq:ref:reference:set:Lambda:introduction} are comparable with each other as limiting tail probabilities evaluated at the rectangles $[0,b]\times[0,1/b]$, $b \in (0,\infty)$, all having the equal volume $b\times (1/b)=1$.
Based on this interpretation, we propose two measures of tail dependence; one is the \emph{maximal tail concordance measure (MTCM)}, defined as the supremum over~\eqref{eq:ref:reference:set:Lambda:introduction}, and another is the $\mu$-\emph{average tail concordance measure (ATCM)}, defined as the normalized average over~\eqref{eq:ref:reference:set:Lambda:introduction} weighted by what we call an \emph{angular measure} $\mu:\mathfrak B((0,\infty))\rightarrow [0,1]$, where $\mathfrak B((0,\infty))$ denotes the Borel $\sigma$-algebra on $(0,\infty)$.
An illustration of the two measures is provided in Figure~\ref{fig:explanations:tail:dependence:measures}.
Inspired by the tail indices proposed by~\citet{furman2015paths}, the MTCM is constructed to extract the most distinctive feature of tail dependence.
On the other hand, the ATCM, regarded as a generalization of the tail dependence measure studied in~\citet{schmid2007multivariateTail}, can be used when specific elements in~\eqref{eq:ref:reference:set:Lambda:introduction} are of particular importance from a practical point of view.
Therefore, the two measures quantify different but important features of non-exchangeable tail dependence summarized by tail copulas.

We also introduce an axiomatic framework of tail dependence measures based on tail copulas as an analog to that introduced in~\citet{scarsini1984} for \emph{measures of concordance}.
The TCMs introduced above are then shown to satisfy all the axiomatic properties naturally required to quantify tail dependence, such as the monotonicity with respect to an appropriate order among tail copulas and the normalization property where the TCM takes the maximum $1$ and the minimum $0$ if and only if the underlying tail copula represents the so-called \emph{tail comonotonicity} and \emph{tail independence}, respectively.
In particular, we investigate the relationship of the proposed measures with tail comonotonicity motivated by~\citet{hua2012taila,hua2012tailb} and~\citet{cheung2019additivity}, the recent studies of tail comonotonicity in the context of risk measures.

Examples of the angular measure $\mu$ and the corresponding ATCM are also provided from the viewpoint of mathematical tractability and practical use.
Our construction of TCMs based on the reference set~\eqref{eq:ref:reference:set:Lambda:introduction} provides a clear probabilistic interpretation of the angular measure $\mu$ and thus enables us to construct flexible tail dependence measures according to the purpose of the analysis.
As an example in credit risk modeling, some elements in~\eqref{eq:ref:reference:set:Lambda:introduction} can be understood as limiting probabilities of joint default when $C$ models the dependence between the values of two firms.
In this situation, the measure $\mu$ can be chosen to put more weight on these elements.

Analytical forms of the proposed measures for a variety of parametric models are
also provided; see Table~\ref{table:analytical:forms}.  Since tail copulas
can be obtained as simple limits of the underlying copula, analytical forms of
the proposed TCMs are typically available.  Admitting analytical forms is
beneficial when parametric copulas are estimated by matching the corresponding
TCMs to their empirical counterparts.
From Table~\ref{table:analytical:forms}, we observe that the MTCM is higher than the TDC for non-exchangeable copulas (1), (4) and (5) whereas the MTCM coincides with the TDC for exchangeable copulas (2) and (3). 
On the other hand, exchangeable copulas (2) and (3) have higher ATCMs compared with others, all of which have the same values of the TDC.

\begin{table}[t!]
  \centering
\resizebox{\columnwidth}{!}{
\begin{tabular}{rccccccl}
\toprule
          &  &  & TDC &  MTCM &  ATCM with $\mu = \frac{1}{2}(\delta_2 + \delta_{\frac{1}{2}})$ \\ [4pt] \midrule
(1)        & \multicolumn{5}{l}{Survival Marshall-Olkin copula $\hat C_{\alpha,\beta}^{\text{MO}}$, $0 < \alpha,\beta \le 1$:}     \\[4pt]
&Analytical form &  &   $\min(\alpha,\beta)$  &   $\sqrt{\alpha\beta}$     &
$ \frac{1}{2}\min(4\alpha,\beta)+\frac{1}{2}\min(\alpha,4\beta)$       \\[3pt]
& $(\alpha,\beta)=(0.35,0.75)$     &  &   0.350  &  {\bf 0.512}  &   0.550  \\ [4pt] \midrule
(2)        & \multicolumn{5}{l}{Archimedean copula $C_{\varphi}$ with $\lim_{x \downarrow 0}x \varphi'(x)/\varphi(x)=-\theta$, $0<\theta< \infty$:}     \\[4pt]
&Analytical form &  &   $2^{-1/\theta}$  &   $2^{-1/\theta}$     &
$ 
2\left(
2^{\theta} + 2^{-\theta}
\right)^{-1/\theta}
$       \\[3pt]
& $\theta = 0.66$     &  &   0.350  &  0.350  &   {\bf 0.600}  \\ [4pt] \midrule
(3)        & \multicolumn{5}{l}{Survival Archimedean copula $\hat C_{\varphi}$ with
$-\lim_{x \downarrow 0}x \varphi'(1-x)/\varphi(1-x)=\theta$, $1\le\theta<\infty$:
}     \\[4pt]
&Analytical form &  &   $2-2^{1/\theta}$  &   $2-2^{1/\theta}$     &
$
  5- 
 2\left(
2^{\theta} + 2^{-\theta}
\right)^{1/\theta}
 $       \\[3pt]
& $\theta = 1.38$     &  &   0.350  &  0.350  &   {\bf 0.584}  \\ [4pt] \midrule
(4)        & \multicolumn{5}{l}{Survival asymmetric Gumbel copula $\hat C_{\alpha,\beta,\theta}^{\text{Gu}}$, $0 < \alpha,\beta\le 1$ and $1\le \theta <\infty$:}     \\[4pt]
&Analytical form &  &   $\alpha + \beta -(\alpha^\theta + \beta^\theta)^{1/\theta}$  & $\left(2-2^{1/\theta}\right)\sqrt{\alpha\beta}$     &
$
\frac{5}{2}(\alpha + \beta)
-
\frac{1}{2}
\{
(4\alpha)^{\theta} + \beta^{\theta}
\}^{1/\theta}
$
      \\[3pt]
      & &  &     &    &
  \hspace{15mm}
  $
-\frac{1}{2}
\{
\alpha^{\theta}
+
(4\beta)^{\theta}
\}^{1/\theta}
$
      \\[3pt]
& $(\alpha,\beta,\theta)=(0.35,0.75,10)$     &&  0.350   & {\bf 0.476}    &   0.550  \\ [4pt] \midrule
(5)        & \multicolumn{5}{l}{Survival asymmetric Galambos copula $\hat C_{\alpha,\beta,\theta}^{\text{Ga}}$, $0 < \alpha,\beta\le 1$ and $0< \theta <\infty$:}     \\[4pt]
&Analytical form &  &   $(\alpha^{-\theta}+\beta^{-\theta})^{-1/\theta}$  &  $2^{-1/\theta}\sqrt{\alpha\beta}$     &
$
\frac{1}{2}
\{
(4\alpha)^{-\theta} + \beta^{-\theta}
\}^{-1/\theta}
$    \\[3pt]
& &  &    &      &
\hspace{15mm}$
+
\frac{1}{2}
\{
\alpha^{-\theta}
+
(4\beta)^{-\theta}
\}^{-1/\theta}
$    \\[3pt]
& $(\alpha,\beta,\theta)=(0.35,0.75,10)$     &  &    0.350 &     {\bf 0.478}  &  0.550   \\ [4pt]
\bottomrule
\end{tabular}
}
\caption{Analytical forms and examples for specific numerical values of
    the TDC and the proposed tail dependence measures; see
    Section~\ref{sec:examples:proposed:tcms} for details.}
\label{table:analytical:forms}
\end{table}

Together with the probabilistic interpretability, we believe that the proposed TCMs can be useful in various practical situations to quantify tail dependence.
Numerical studies are also conducted to demonstrate the practical use of the proposed TCMs.
Our simulation studies show that different features of tail dependence are captured by the proposed two measures.
A real data analysis then reveals striking tail dependence and tail non-exchangeability of the return series of stock indices particularly in the periods of financial distress.

The present paper is organized as follows.
Section~\ref{sec:tail:copulas:tail:concordance:measures} presents the framework for measuring tail dependence based on tail copulas.
Section~\ref{sec:proposed:tcm} introduces the proposed tail concordance measures.
Axiomatic properties and examples for various parametric copulas are also provided.
Simulation and empirical studies are conducted in Section~\ref{sec:numerical:experiments}.
Section~\ref{sec:conclusion} concludes with ideas for future research directions.
Technical supplements, proofs and statistical inference of the proposed measures can be found in the appendix.

\section{Tail copulas and tail concordance measures}\label{sec:tail:copulas:tail:concordance:measures}

Let $\mathcal C_2$ be the set of all $2$-\emph{copulas}, that is, bivariate distribution functions with standard uniform margins.
The \emph{comonotonicity} and \emph{independence copulas} are defined by $M(u,v)=\min(u,v)$ and $\Pi(u,v)=uv$ for $(u,v)\in [0,1]^2$, respectively.
For $(U,V)\sim C \in \mathcal C_2$, the limit
\begin{align}\label{eq:def:tail:copulas}
\Lambda(u,v)=\Lambda(u,v;C)
=\lim_{p\downarrow 0}\frac{C(pu,pv)}{p}
=\lim_{p\downarrow 0}\frac{\Prob((U,V)/p \in [0,u]\times[0,v])}{p},\quad (u,v) \in \IR_{+}^2,
\end{align}
provided it exists, is called a \emph{tail copula}, also known as a \emph{tail dependence function}.
Throughout this paper, we focus only on the lower tail around the origin $(0,0)$ in $[0,1]^2$ since tail dependence around the other three corners $(1,0)$, $(0,1)$, $(1,1)$ can be studied by replacing $C$ with its rotated copulas.
The existence of $\Lambda(\cdot;C)$ can equivalently be stated in terms of the so-called \emph{tail expansion} of $C$; see \cite{jaworski2004uniform,jaworski2006uniform,jaworski2010tail} for details.
Basic properties of tail copulas are summarized in Appendix~\ref{sec:basic:properties:tail:copulas}.
Let
\begin{align*}
\mathcal C_2^{\text{L}}=\left\{C\in\mathcal C_2: \lim_{p\downarrow0}\frac{C(pu,pv)}{p}\text{ exists for every } (u,v) \in \IR_{+}^2\right\}
\end{align*}
be the set of all copulas admitting tail copulas, and let
\begin{align*}
\mathcal L=\left\{\Lambda:\IR_{+}^2\rightarrow \IR_{+}: \text{ there exists $C\in \mathcal C_2$ such that }\Lambda(u,v)=\lim_{p\downarrow 0}\frac{C(pu,pv)}{p}\right\}
\end{align*}
be the set of all tail copulas.
By construction, $\mathcal C_2^{\text{L}}$ and $\mathcal L$ are convex sets.
Moreover, the inclusion relationship $\mathcal C_2^{\text{L}}\subseteq \mathcal C_2$ is strict; see \citet[Corollary~8.3.2]{jaworski2010tail}.

The following concepts are fundamental for quantifying tail dependence summarized by tail copulas.
\begin{definition}[Tail dependence and tail concordance order]\label{def:tail:dependence:tail:concordance:order}
Let $C,\,C' \in \mathcal C_2^\text{L}$ with $\Lambda=\Lambda(\cdot;C)$ and $\Lambda'=\Lambda(\cdot;C')$.
\begin{enumerate}[label=\arabic*),itemsep=0pt]
\item\label{item:def:tail:independence} {\bf (Tail independence)}
$C$ (or $\Lambda$) is called \emph{tail independent} if $\Lambda\equiv 0$, and it is called \emph{tail dependent} if $\Lambda \not\equiv 0$.
\item\label{item:def:tail:comonotonic} {\bf (Tail comonotonicity)}
$C$ (or $\Lambda$) is called \emph{tail comonotonic} if $\Lambda=\overline \Lambda$, where
\begin{align*}
\overline \Lambda(u,v)=\min(u,v),\quad (u,v)\in \IR_{+}^2.
\end{align*}
\item\label{item:def:tail:concordance:order} {\bf (Tail concordance order)} $C'$ (or $\Lambda'$) is said to be more \emph{tail concordant} than $C$ (or $\Lambda$), denoted by $C \preceq_\text{L} C'$ (or $\Lambda \preceq \Lambda'$), if $\Lambda(u,v)\leq \Lambda'(u,v)$ for all $(u,v)\in \IR_{+}^2$.
\end{enumerate}
\end{definition}
Note that $\Pi$ is tail independent and $M$ is tail comonotonic.
Therefore, we have that $0, \overline \Lambda \in \mathcal L$.
Roughly speaking, tail comonotonicity means that the underlying copula behaves like the comonotonic copula $M$ in the tail.
Detailed discussions on the notion of tail comonotonicity, particularly related to asymptotic behavior of risk measures, can be found in~\citet{hua2012taila,hua2012tailb} and~\citet{cheung2019additivity}.
On the other hand, tail independence means that the joint probability $C(pu,pv)$ vanishes faster than $p$.
By Proposition~\ref{prop:basic:properties:tail:dependence:function}~\ref{basic:prop:item:degeneracy}, a tail copula $\Lambda$ is tail independent if and only if there exists $u_0,v_0>0$ such that $\Lambda(u_0,v_0)=0$.
Moreover, as stated in Proposition~\ref{prop:basic:properties:tail:dependence:function}~\ref{basic:prop:item:bounds}, the tail copulas $0$ and $\overline \Lambda$ are the minimal and maximal elements in $\mathcal L$ with respect to the tail concordance order, that is, $0\preceq \Lambda \preceq \overline \Lambda$ for all $\Lambda \in \mathcal L$.
In fact, \cite{jaworski2004uniform} showed that $\mathcal L$ can be characterized by
\begin{align*}
\mathcal L=\{\Lambda:\IR_{+}^2\rightarrow \IR_{+}: \Lambda \text{ is positive homogeneous, $2$-increasing and } 0\preceq \Lambda\preceq \overline \Lambda\};
\end{align*}
see Proposition~\ref{prop:basic:properties:tail:dependence:function}~\ref{basic:prop:item:2:increasing} and~\ref{basic:prop:item:one:homogeneity} for the properties of 2-increasingness and positive homogeneity of tail copulas.
A key property of tail copulas is the positive homogeneity
\begin{align}\label{eq:ref:1:homogeneity}
\Lambda(tu,tv)=t\Lambda(u,v)\quad\text{for every}\quad t \geq 0 \quad\text{and}\quad (u,v)\in \IR_{+}^2,
\end{align}
which leads to the following equivalent relations to the tail concordance order.
\begin{proposition}[Properties related to the tail concordance order]
\label{prop:tail:concordance:order:tail:probabilities}
Let $C,\,C' \in \mathcal C_2^\text{L}$ with $\Lambda=\Lambda(\cdot;C)$ and $\Lambda'=\Lambda(\cdot;C')$.
\begin{enumerate}[label=\arabic*), itemsep=0pt]
\item\label{item:tail:concordance:order:theta:angle}
Let $\theta \mapsto r_\theta \in (0,\infty)$, $\theta \in (0,\pi/2)$, be arbitrary.
Then $\Lambda \preceq\Lambda'$ if and only if
\begin{align}\label{eq:simplified:condition:Lambda:concordance:order:theta:angle}
\Lambda(r_\theta\cos\theta,r_\theta\sin\theta)\leq \Lambda'(r_\theta\cos\theta,r_\theta\sin\theta)\quad\text{for every} \quad \theta \in \left(0,\frac{\pi}{2}\right).
\end{align}
\item\label{item:tail:concordance:order:tail:probability}
We have that $C \preceq_\text{L} C'$ if and only if, for $(U,V)\sim C$ and $(U',V')\sim C'$,
\begin{align*}
\lim_{p\downarrow 0}\Prob(U\leq pu\,|\, V \leq p)&\leq \lim_{p\downarrow 0}\Prob(U'\leq pu\,|\, V' \leq p)\quad\text{for every}\quad u \in (0,1],\quad\text{and}\\
\lim_{p\downarrow 0}\Prob(V\leq pv\,|\, U \leq p)&\leq \lim_{p\downarrow 0}\Prob(V'\leq pv\,|\, U' \leq p)\quad\text{for every}\quad v \in (0,1].
\end{align*}
\end{enumerate}
\end{proposition}
By Proposition~\ref{prop:tail:concordance:order:tail:probabilities}~\ref{item:tail:concordance:order:theta:angle}, the relationship $\Lambda\preceq \Lambda'$ can be simplified to checking the pointwise inequality between $\Lambda$ and $\Lambda'$ at one representative point $(r_\theta\cos\theta,r_\theta\sin\theta)$ for each angle $\theta \in (0,\pi/2)$.
Part~\ref{item:tail:concordance:order:tail:probability} provides an intuitive interpretation of the order $\Lambda \preceq \Lambda'$ by the monotonicity of limiting tail probabilities.
Related orders implied by $\Lambda\preceq \Lambda'$ can be found in \cite{li2013dependence}.

\begin{remark}
Independently of the present paper, the tail concordance order has recently been introduced and explored in~\citet{siburg2022comparing} and~\citet{siburg2022multivariate}.
\end{remark}

\begin{example}[Tail concordance order for extreme value copulas]
\label{example:tail:concordance:order:ev:copulas}
A bivariate \emph{extreme value (EV) copula} is given by
\begin{align*}
C_A(u,v)=\exp\left\{(\log u + \log v)A\left(\frac{\log u}{\log u + \log v}\right)\right\},\quad (u,v)\in [0,1]^2,
\end{align*}
where $A \in \mathcal A$ is the so-called \emph{Pickands dependence function} with
\begin{align*}
  \mathcal A = \{A:[0,1]\rightarrow [1/2,1]: \text{convex and }\max(w,1-w)\leq A(w)\leq 1\text{ for all }w \in [0,1]\}.
\end{align*}
For instance, $A\equiv 1$ yields $C_A=\Pi$ and $A(w)=\max(w,1-w)$ yields $C_A=M$.
Using the relationships $\log(1-x)\simeq -x$ and $1-x \simeq e^{-x}$ ($x \downarrow 0$), it is straightforward to check that the survival copula of $C_A$ has the tail copula
\begin{align}\label{eq:ev:copula:tail:dependence:function}
\Lambda(u,v;\hat C_A)=u+v-(u+v)A\left(\frac{u}{u+v}\right),\quad (u,v)\in\IR_{+}^2.
\end{align}
By Equation~\eqref{eq:ev:copula:tail:dependence:function}, the tail concordance order $C_A\preceq_\text{L} C_{A'}$ for $A, A' \in \mathcal A$ is equivalent to $A'(w)\leq A(w)$ for all $w \in [0,1]$, which is the order considered in~\cite{jaworski2019extreme} to quantify dependence of EV copulas.
\end{example}

We now introduce axioms of measures that quantify the degree of tail concordance.
\begin{definition}[Axioms of tail concordance measures]
\label{def:axioms:tail:concordance:measure}
A map $\kappa:\mathcal L \rightarrow \IR$ is called a \emph{tail concordance measure (TCM)} if it satisfies the following conditions.
\begin{enumerate}[label=\arabic*), itemsep=0pt]
    \item\label{def:axioms:item:normalization} {\bf (Normalization)} $\kappa(\Lambda)=1$ if $\Lambda$ is tail comonotonic.
    \item\label{def:axioms:item:tail:independence} {\bf (Tail independence)} $\kappa(\Lambda)=0$ if and only if $\Lambda$ is tail independent.
    \item\label{def:axioms:item:monotonicity} {\bf (Monotonicity)} If $\Lambda\preceq \Lambda'$ for $\Lambda,\Lambda' \in \mathcal L$, then $\kappa(\Lambda)\leq \kappa(\Lambda')$.
    \item\label{def:axioms:item:continuity} {\bf (Continuity)}
    If $\Lambda^{[n]}\to\Lambda$ ($n\to\infty$) pointwise for $\Lambda^{[n]},\Lambda \in \mathcal L$, $n\in\IN$,
    then
    $\lim_{n\rightarrow \infty}\kappa(\Lambda^{[n]})=\kappa(\Lambda)$.
\end{enumerate}
If, in addition, $\kappa(\Lambda)=1$ and $\Lambda=\overline \Lambda$ are equivalent, then $\kappa$ is called a \emph{strict} TCM.
\end{definition}
Axiom~\ref{def:axioms:item:monotonicity} is a fundamental requirement to quantify tail concordance, and Axioms~\ref{def:axioms:item:normalization} and~\ref{def:axioms:item:tail:independence} normalize the measure.
By these three axioms, we have that $0\leq \kappa(\Lambda)\leq 1$ for every $\Lambda \in \mathcal L$.
Axiom~\ref{def:axioms:item:continuity} ensures that if $\kappa$ is calculated based on an approximated tail copula $\Lambda^{[n]}$ of $\Lambda$, then the estimate $\kappa(\Lambda^{[n]})$ is close to $\kappa(\Lambda)$.
Finally, strictness of $\kappa$ can be naturally required to detect tail comonotonicity under which various asymptotic results are accessible.
Let us now introduce properties of TCMs related to the evaluation of convex combinations of tail copulas.

\begin{definition}[Convexity, concavity and linearity]
\label{def:tcm:strictness:convexity:concavity}
A TCM $\kappa$ is called
\begin{enumerate}[label=\arabic*), itemsep=0pt]
\item\label{item:properties:convex} \emph{convex} if
\begin{align}\label{eq:def:convexity:xi}
\kappa(t\Lambda+(1-t)\Lambda')\leq t\kappa(\Lambda)+(1-t)\kappa(\Lambda')\quad\text{for every}\quad \Lambda,\Lambda' \in \mathcal L\quad\text{and}\quad t \in [0,1],
\end{align}
\item\label{item:properties:concave} \emph{concave} if the reverse inequality in~\eqref{eq:def:convexity:xi} holds, and
\item\label{item:properties:linear} \emph{linear} if $\kappa$ is convex and concave.
\end{enumerate}
\end{definition}

\section{The proposed tail concordance measures}
\label{sec:proposed:tcm}

\subsection{Definitions, basic properties and examples}\label{sec:definition:proposed:tcms:examples}

By positive homogeneity~\eqref{eq:ref:1:homogeneity} of tail copulas, it is sufficient to construct measures of tail concordance only from $\Lambda(u,v)$ with the restricted domain $(u,v) \in \{(r_\theta\cos\theta,r_\theta\sin\theta): \theta \in (0,\pi/2)\}$ for some $r_\theta:(0,\pi/2)\rightarrow (0,\infty)$.
As discussed in Section~\ref{sec:introduction}, we take $r_\theta = 1/(\sin\theta \cos \theta)^{1/2}$ so that the rectangles appearing in~\eqref{eq:def:tail:copulas} have the same volume as $[0,1]\times[0,1]$ for every $\theta \in (0,\pi/2)$.
Consequently, elements in the reference set
\begin{align*}
\left\{\Lambda\left(b,\frac{1}{b}\right): b \in (0,\infty)\right\}
\end{align*}
are comparable with each other as limiting tail probabilities evaluated at the rectangles $[0,b]\times[0,1/b]$ for $b \in (0,\infty)$.
Based on this interpretation, we consider the following classes of measures.

\begin{definition}[Maximal and $\mu$-average tail concordance measures]
\label{def:maximal:mu:average:tcms}
\hspace{2mm}
\begin{enumerate}[label=\arabic*), itemsep=0pt]
\item The \emph{maximal TCM (MTCM)} $\lambda^\ast:\mathcal L \rightarrow \IR_{+}$ is defined by
\begin{align}\label{eq:def:mtcm}
\lambda^\ast(\Lambda)=\sup_{b \in (0,\infty)}\Lambda(b,1/b).
\end{align}
\item Let $\mathcal M$ be the set of all Borel probability measures on $(0,\infty)$ such that $\int_{(0,\infty)}\overline \Lambda(b,1/b)\rd \mu(b)<\infty$.
For $\mu \in \mathcal M$, the $\mu$-\emph{average TCM (ATCM)} $\lambda_\mu:\mathcal L \rightarrow \IR_{+}$ is defined by
\begin{align}\label{eq:def:mu:atcm}
\lambda_\mu(\Lambda)=\frac{\int_{(0,\infty)} \Lambda(b,1/b)\rd \mu(b)}{\int_{(0,\infty)}\overline \Lambda(b,1/b)\rd \mu(b)}.
\end{align}
The measure $\mu$ is called the \emph{angular measure}.
\end{enumerate}
\end{definition}

The MTCM is constructed based on the idea that the rectangle maximizing the limiting probability $\Lambda(b,1/b)$ captures the feature of tail dependence of $C$.
On the other hand, the ATCM can be interpreted as a normalized average of the limiting tail probabilities weighted by $\mu \in \mathcal M$.
As we will see in Example~\ref{example:mu:atcms:practical}, the angular measure $\mu \in \mathcal M$ can be chosen externally by an analyst depending on the importance of each rectangle $[0,b]\times[0,1/b]$.
Therefore, the two measures $\lambda^\ast$ and $\lambda_\mu$ can be used for different purposes since the main objective of $\lambda^\ast$ may be to extract the most distinctive feature of the tail behavior.
Unlike the MTCM, the ATCM requires the denominator $\int_{(0,\infty)}\overline \Lambda(b,1/b)\rd \mu(b)$ so that $\lambda_\mu$ satisfies Axiom~\ref{def:axioms:item:normalization} in Definition~\ref{def:axioms:tail:concordance:measure}; see Proposition~\ref{prop:axiomatic:properties:proposed:tcms} below.
Due to this denominator, it may be more reasonable to use the numerator in~\eqref{eq:def:mu:atcm} rather than $\lambda_\mu$
when quantifying the degree of average tail dependence in comparison with $\lambda$ and $\lambda^\ast$.
Differently from the TDC, $\mu$-ATCMs can capture non-exchangeable tail dependence as seen in the following representation.

\begin{proposition}[Representation of $\mu$-ATCMs]
\label{prop:representation:mu:atcm}
For $\mu \in \mathcal M$, the $\mu$-ATCM can be represented as
\begin{align}\label{eq:mu:atcm:representation}
\lambda_\mu(\Lambda)=
 \frac{w_0 \Lambda(1,1)+w_1 \int_{(0,1)}\Lambda(b,1/b)\rd \mu_1(b)+w_2\int_{(0,1)}\Lambda(1/b,b)\rd\mu_2(b)}{w_0 +w_1 \int_{(0,1)}b\rd \mu_1(b)+w_2\int_{(0,1)}b\rd \mu_2(b)},
\end{align}
where $w_0,w_1,w_2\geq 0$ with $w_0+w_1+w_2=1$ and $\mu_1$ and $\mu_2$ are Borel probability measures on $(0,1)$.
\end{proposition}

Representation~\eqref{eq:mu:atcm:representation} may be easier to interpret than~\eqref{eq:def:mu:atcm} since the limiting non-exchangeable tail probabilities $\Lambda(b,1/b)$ and $\Lambda(1/b,b)$ are treated on the same scale on $(0,1)$ where the weights for different $b$s are determined by $\mu_1$ and $\mu_2$, respectively.
The TDC corresponds to the case when $(w_0,w_1,w_2)=(1,0,0)$, that is, the weight is concentrated on the main diagonal.
By taking positive values for $w_1$ and $w_2$, the $\mu$-ATCM incorporates the limiting tail probabilities $\Lambda(b,1/b)$ and $\Lambda(1/b,b)$, which are in general different for non-exchangeable copulas.

\begin{remark}[Attainability of the MTCM]\label{remark:attainability:mtcm}
The supremum in~\eqref{eq:def:mtcm} is attainable in $b \in (0,\infty)$ for any $\Lambda \in \mathcal L$ since $\lim_{b\downarrow 0}\Lambda\left(b,1/b\right)=\lim_{b\rightarrow \infty}\Lambda\left(b,1/b\right)=0$ and the map $b \mapsto \Lambda\left(b,1/b\right)$ is continuous and bounded; see Proposition~\ref{prop:basic:properties:tail:dependence:function}~\ref{basic:prop:item:bounds} and Proposition~\ref{prop:continuity:derivatives}~\ref{basic:prop:item:continuity}.
Therefore, one can write $\lambda^\ast(\Lambda)=\max_{b \in (0,\infty)}\Lambda(b,1/b)$.
Moreover, assuming that the maximum is uniquely attained at a single point, we write
\begin{align*}
b^\ast=\argmax_{b \in (0,\infty)}\Lambda(b,1/b).
\end{align*}
A deviation of $b^\ast$ from $1$ may be an important sign of tail non-exchangeability under which the TDC may not be a suitable measure to summarize tail dependence; see Section~\ref{subsec:real:data:analysis} for experiments.
\end{remark}

\begin{remark}
Independently of the present paper,~\citet{siburg2022comparing} introduced tail dependence measures constructed as an
average and maximum over the reference set $\{\Lambda(s,1-s): s \in [0,1] \}$.
\end{remark}

As seen in the following examples, the angular measure $\mu$ of the $\mu$-ATCM can be chosen for mathematical tractability or for practical purposes.
\begin{example}[Examples of $\mu$-ATCMs]
\label{example:mu:atcms}
\hspace{2mm}
\begin{enumerate}[label=\arabic*), itemsep=0pt]
\item\label{example:item:generalized:tdc} \emph{Generalized TDC}:
For $b \in (0,\infty)$, let $\mu=\delta_b$ in~\eqref{eq:def:mu:atcm} be the Dirac measure on $b$.
Then the resulting $\mu$-ATCM is given by
    \begin{align*}
    \lambda_{\delta_b}(\Lambda)=\frac{\Lambda(b,1/b)}{\overline \Lambda(b,1/b)}=
    \Lambda\left(1\vee b^2,1\vee\frac{1}{b^2}\right).
    \end{align*}
    We call this measure the \emph{generalized tail dependence coefficient (GTDC)}; note that $b=1$ leads to the TDC.
\item\label{example:item:uniform:atcm}  \emph{Uniform ATCM}: Let $w_1=w_2=1/2$, and $\mu_1=\mu_2$ in~\eqref{eq:mu:atcm:representation} be probability measures whose unit masses are uniformly put on $(0,1)$.
Then the corresponding $\mu$-ATCM, denoted by $\lambda_\text{U}$, is given by
\begin{align*}
\lambda_\text{U}(\Lambda)=\int_{(0,1)}\left\{\Lambda\left(b,\frac{1}{b}\right)+\Lambda\left(\frac{1}{b},b\right) \right\}\rd b.
\end{align*}
We call this measure the \emph{uniform ATCM}.
\end{enumerate}
\end{example}

\begin{example}[A practical choice of $\mu$ in credit risk modeling]
\label{example:mu:atcms:practical}
Let $X \sim F$ and $Y \sim G$ be continuously distributed random variables representing the values of two firms with default probability of the former being estimated as $0.05$ and that of the latter being estimated as an interval $[0.001, 0.01]$ with some levels of credibility.
To model the dependence between $X$ and $Y$, joint default events of the form
\begin{align*}
\{X\leq q_{0.05}(X),\,Y \leq q_\alpha(Y)\}=
\{U\leq 0.05,\,V \leq \alpha\},\quad \alpha \in [0.001,\, 0.01],
\end{align*}
are of primary concern, where $q_\alpha(X)$ and $q_\alpha(Y)$, $\alpha \in (0,1)$, are the $\alpha$-quantiles of $F$ and $G$, and $(U,V)=(F(X),G(Y))$.
Therefore, to quantify the tail dependence of $(X,Y)$, the values of the tail copulas
\begin{align*}
\Lambda(b,1/b),\quad b \in \left[\left(\frac{0.001}{0.05}\right)^{1/2},\left(\frac{0.01}{0.05}\right)^{1/2}\right]=[0.141,0.447],
\end{align*}
may be more important than $\Lambda(1,1)$, where $b \in [0.141,0.447]$ is determined so that the ratio between $0.05$ and $\alpha \in [0.001,0.01]$ equals that of $pb$ and $p/b$.
Therefore, compared with the TDC, the tail dependence of interest may be better summarized by a $\mu$-ATCM with $\mu$ supported on $[0.141,0.447]$, and with the weights possibly determined proportionally to the credibility of the estimated default probabilities.
\end{example}

The next proposition states that the MTCM and ATCMs are TCMs in the sense of Definition~\ref{def:axioms:tail:concordance:measure}.

\begin{proposition}[Axiomatic properties]
\label{prop:axiomatic:properties:proposed:tcms}
\hspace{2mm}
\begin{enumerate}[label=\arabic*), itemsep=0pt]
\item\label{prop:item:mtcm:is:convex:tcm} The MTCM $\lambda^\ast$ is a strict and convex TCM.
\item\label{prop:item:atcm:is:linear:tcm} For $\mu \in \mathcal M$, the $\mu$-ATCM $\lambda_\mu$ is a linear TCM.
\item\label{prop:item:condition:atcm:is:strict:tcm} A $\mu$-ATCM is strict if $\mu \in \mathcal M$ satisfies the following condition:
    \begin{align}\label{eq:condition:probability:around:diagonal}
    \mu((1-\epsilon,1+\epsilon)\cap(0,\infty))>0\text{ for any } \epsilon >0.
    \end{align}
\end{enumerate}
\end{proposition}

Proposition~\ref{prop:axiomatic:properties:proposed:tcms}~\ref{prop:item:mtcm:is:convex:tcm} and~\ref{prop:item:atcm:is:linear:tcm} state that the proposed measures satisfy the axiomatic properties of TCMs presented in Definition~\ref{def:axioms:tail:concordance:measure}.
According to Proposition~\ref{prop:axiomatic:properties:proposed:tcms}~\ref{prop:item:condition:atcm:is:strict:tcm}, the uniform ATCM $\lambda_{\text{U}}$ in Example~\ref{example:mu:atcms}~\ref{example:item:uniform:atcm} and the TDC $\lambda$ are strict ATCMs.
However, as seen in the following example, the GTDC in Example~\ref{example:mu:atcms}~\ref{example:item:generalized:tdc} is not always a strict TCM as the corresponding angular measure $\mu$ may violate Condition~\eqref{eq:condition:probability:around:diagonal}.

\begin{example}[GTDCs of a singular copula]
\label{example:atcm:nelsen:singular:copula}
  For $\theta \in (0,1)$, let $C_\theta$ be a copula considered in \citet[Section~3.2.1]{nelsen2006introduction}, where the
  probability mass $\theta \in (0,1)$ is uniformly distributed on the line segment from
  $(0,0)$ to $(\theta,1)$, and the probability mass $1-\theta$ is uniformly
  distributed on the line segment from $(\theta,1)$ to $(1,0)$.
  The tail copula of $C_\theta$ is given by
    $\Lambda_\theta(u,v)=\Lambda(u,v;C_\theta)=\min(u,\theta v)$,
  and thus the GTDC for $\mu=\delta_b$ is given by
  \begin{align*}
    \lambda_{\delta_b}(\Lambda_\theta)=
    \begin{cases}
      1, & \text{ if } b\leq\sqrt{\theta},\\
      \theta/b^2 \in (\theta,1), & \text{ if } \sqrt{\theta}<b< 1, \\
      \theta, & \text{ if } 1\leq b.\\
    \end{cases}
  \end{align*}
  Therefore, $\Lambda_\theta \neq \overline \Lambda$ but it holds that $\lambda_{\delta_{\sqrt{\theta}}}(\Lambda_\theta)=1$.
\end{example}

\subsection{Bounds of average tail concordance measures}\label{sec:minimal:maximal:average:tail:dependence:measure}

We next study bounds of $\mu$-ATCMs for a fixed $\Lambda \in \mathcal L$ over all angular measures $\mu \in \mathcal M$.
Namely, we are interested in the \emph{minimal} and \emph{maximal average tail concordance measures} defined by
\begin{align*}
\underline \lambda(\Lambda)=\inf_{\mu\in \mathcal M}\lambda_\mu(\Lambda)
\quad\text{and}\quad
\overline \lambda(\Lambda)=\sup_{\mu\in \mathcal M}\lambda_\mu(\Lambda)
,\quad \Lambda \in \mathcal L,
\end{align*}
and their attaining angular measures.
Characterizations of these extremal ATCMs and their connections to the MTCM are provided in the following theorem.

\begin{theorem}[Characterization of the minimal and maximal ATCMs]
\label{thm:characterization:minimal:maximal:atcm}
Let $\Lambda \in \mathcal L$.
\begin{enumerate}[label=\arabic*), itemsep=0pt]
\item \label{thm:item:characterization:minimal:TCM}
For every $\mu \in \mathcal M$, it holds that
\begin{align}\label{eq:inequality:tail:dependence:coefficient}
\lambda(\Lambda)\leq \lambda_{\mu}(\Lambda)\leq 1\wedge a\lambda(\Lambda)\quad\text{ where }\quad
    a=a(\mu)=\frac{\int_{(0,\infty)}\max(b,1/b)\rd \mu(b)}{\int_{(0,\infty)}\overline \Lambda(b,1/b)\rd \mu(b)} \in [1,\infty].
\end{align}
Therefore, the angular measure $\delta_1$ attains $\underline \lambda(\Lambda)=\lambda(\Lambda)=\Lambda(1,1)$.
\item\label{thm:item:characterization:maximal:average:tail:concordance:measure}
It holds that
\begin{align}\label{eq:formula:maximal:atcm}
\overline \lambda(\Lambda)&=\sup_{b \in (0,\infty)}\frac{
\Lambda\left(b,1/b\right)
}{
\overline \Lambda\left(b,1/b\right)
}\\
\label{eq:formula:maximal:atcm:max}&=
\max\left(\lim_{b \downarrow 0}\Lambda\left(1,\frac{1}{b^2}\right),\ \lim_{b \rightarrow \infty}\Lambda\left(b^2,1\right)\right).
\end{align}
The supremum in~\eqref{eq:formula:maximal:atcm} cannot be replaced by the maximum in general.
\item\label{thm:item:bounds:mtcm}
It holds that $\lambda(\Lambda)\leq \lambda^\ast(\Lambda)\leq \overline \lambda(\Lambda)$.
\end{enumerate}
\end{theorem}

Theorem~\ref{thm:characterization:minimal:maximal:atcm}~\ref{thm:item:characterization:minimal:TCM} shows that the TDC is the minimal ATCM.
By~\eqref{eq:inequality:tail:dependence:coefficient}, the coefficient
$a=a(\mu)$ quantifies the gap between $\lambda$ and $\lambda_\mu$.
For the GTDC in Example~\ref{example:mu:atcms}~\ref{example:item:generalized:tdc}, we have that $a(\delta_b)=\max(b^2,1/b^2)$.
Therefore, $a(\delta_1)=1$, and $a(\delta_b)$ tents to inifinity as $b \downarrow 0$ or $b \rightarrow \infty$.
Moreover, $a$ may admit the value $\infty$ when, for example, $\mu$ corresponds to the uniform ATCM in Example~\ref{example:mu:atcms}~\ref{example:item:uniform:atcm}.

Theorem~\ref{thm:characterization:minimal:maximal:atcm}~\ref{thm:item:characterization:maximal:average:tail:concordance:measure} provides formulas for the maximal ATCM.
By Formula~\eqref{eq:formula:maximal:atcm}, the maximal ATCM can be interpreted as finding a rectangle maximizing the limiting tail probability $\Lambda(b,1/b)$ normalized by $\overline \Lambda(b,1/b)$.
Despite this intuitive interpretation, the maximal ATCM is not an appealing measure of tail dependence as we will now explain.
First and foremost, maximizing the ratio $\Lambda(b,1/b)/\overline \Lambda(b,1/b)$ does not extract informative features of the underlying tail dependence since $b\mapsto \Lambda(b,1/b)/\overline \Lambda(b,1/b)$ is decreasing on $(0,1)$ and increasing on $(1,\infty)$; therefore, as seen in Formula~\eqref{eq:formula:maximal:atcm:max}, the maximal ATCM can be interpreted as being attained at the $y$-axis ($b\downarrow 0$) or the $x$-axis ($b \rightarrow \infty$) regardless of $\Lambda$.
In fact, the maximal ATCM is typically independent of the parameters of tail dependence of various parametric copulas.
For example, the maximal ATCMs of the survival Gumbel and Clayton copulas are all $1$ regardless of their parameters; see Section~\ref{sec:examples:proposed:tcms}.
In summary, although the maximal ATCM provides some insights on the choice of $\mu$, it is not recommendable to use as a tail dependence measure.

Finally, Theorem~\ref{thm:characterization:minimal:maximal:atcm}~\ref{thm:item:bounds:mtcm} says that the MTCM $\lambda^\ast$ is also bounded by $\lambda$ and $\overline \lambda$ although the MTCM is not an ATCM.

In the next remark we shall adopt the notation $\mathrm{D}_j f(\bx)=\partial f(\bx)/\partial x_j$ for $\bx=(\sqc{x}{d})$ and $f:\IR^d\rightarrow \IR$, provided the partial derivative exists.

\begin{remark}[Attainability of the maximal ATCM]\label{remark:attainability]maximal:atcm}
As seen in Proposition~\ref{prop:basic:properties:tail:dependence:function}~\ref{basic:prop:item:monotonicity} and Proposition~\ref{prop:continuity:derivatives}~\ref{basic:prop:item:partial:derivatives}, the maps $t \mapsto \Lambda(t,1)$ and $t\mapsto \Lambda(1,t)$, $t \in \IR_{+}$, are increasing.
Therefore, by Theorem~\ref{thm:characterization:minimal:maximal:atcm}~\ref{thm:item:characterization:maximal:average:tail:concordance:measure}, the maximal ATCM is attained at $\overline \mu=\delta_b \in \mathcal M$ for some $b \in (0,\infty)$ if and only if there exists $t_0>0$ such that
$\mathrm{D}_1\Lambda(t,1)=0$ for $t \geq t_0$ when $\overline \lambda(\Lambda)=\lim_{t\rightarrow \infty}\Lambda(t,1)$, and $\mathrm{D}_2\Lambda(1,t)=0$ for $t \geq t_0$ when $\overline \lambda(\Lambda)=\lim_{t\rightarrow \infty}\Lambda(1,t)$.
An example of such an attainable case can be found in Example~\ref{example:marshall:olkin:copula:maximal:tail:dependence:measure}.
\end{remark}

\subsection{Illustrative examples}\label{sec:examples:proposed:tcms}

In this section we derive the proposed TCMs for various parametric copulas.

\begin{example}[Survival Marshall--Olkin copula]\label{example:marshall:olkin:copula:maximal:tail:dependence:measure}
Let us consider the \emph{Marshall--Olkin copula} $C_{\alpha,\beta}$ defined by
\begin{align*}
C_{\alpha,\beta}^{\text{MO}}(u,v)=\min\left(u^{1-\alpha}v, uv^{1-\beta}\right),\quad
\alpha,\beta \in (0,1],\quad (u,v)\in [0,1].
\end{align*}
By calculation, the tail copula of the survival Marshall--Olkin copula is given by
\begin{align*}
\Lambda_{\alpha,\beta}(u,v)=\Lambda(u,v;\hat C_{\alpha,\beta})=u+v-\max(v+(1-\alpha)u,u+(1-\beta)v)=\min(\alpha u,\beta v).
\end{align*}
Therefore, the TDC is given by $\lambda(\Lambda_{\alpha,\beta})=\Lambda_{\alpha,\beta}(1,1)=\min(\alpha,\beta)$.
Moreover, since the function $b\mapsto \Lambda\left(b,1/b\right)=\min\left(\alpha b, \beta/b\right)$, $b \in (0,\infty)$, is maximized at $b^\ast=\sqrt{\beta/\alpha}$, we have that
$\lambda^\ast(\Lambda_{\alpha,\beta})=\sqrt{\alpha\beta}$.
Finally, by~\eqref{eq:formula:maximal:atcm:max}, we have that $\overline \lambda(\Lambda_{\alpha,\beta})=\max(\alpha,\beta)$, which is attainable at $\overline \mu=\delta_b$ for every $b \in \left(0,\sqrt{\beta/\alpha}\right]$ if $\alpha \geq \beta$, and for every $b \in \left[\sqrt{\beta/\alpha},\infty\right)$ if $\alpha \leq \beta$.
\end{example}

\begin{figure}[t!]
  \centering
  \includegraphics[width=16 cm]{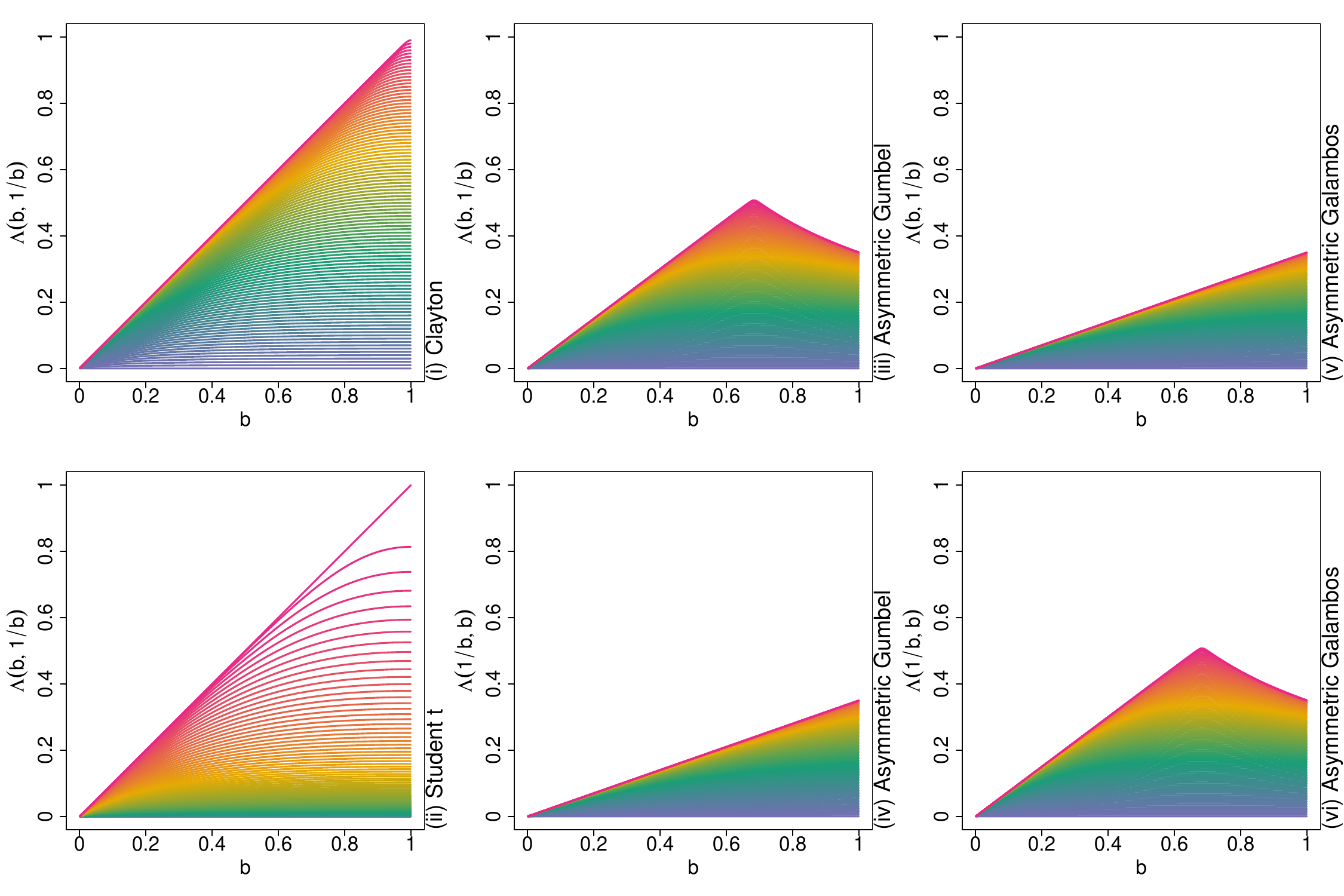}
 \caption{
 The tail copulas $b \mapsto \Lambda\left(b,1/b\right)$ and $b \mapsto \Lambda\left(1/b,b\right)$ for a (i) Clayton copula $C_\theta^{\text{Cl}}$; (ii) $t$ copula $C_{\nu,\rho}^t$ with $\nu=5$~\citep{nikoloulopoulos2009extreme}; (iii), (iv) survival asymmetric Gumbel copula $C_{\alpha,\beta,\theta}^\text{Gu}$ with $\alpha=0.75$ and $\beta=0.35$; and (v), (vi) survival asymmetric Galambos copula $C_{\alpha,\beta,\theta}^\text{Ga}$ with $\alpha=0.35$ and $\beta=0.75$.
The parameters of the copulas vary from (i) $\theta=0.075$ to $68.967$; (ii) $\rho=-0.99$ to $1$; (iii), (iv) $\theta=1.00$ to $69.661$; and (v), (vi) $\theta=0.075$ to $68.967$ as colors vary from blue to green to yellow, and to red.
For Clayton and $t$ copulas, only $b \mapsto \Lambda\left(b,1/b\right)$, $b \in [0,1]$, is plotted since these copulas are exchangeable.
  }
  \label{fig:example:tail:dependence:functions}
\end{figure}

\begin{example}[Archimedean copulas]\label{example:archimedean:copulas}
Consider an \emph{Archimedean copula}
\begin{align*}
C_\varphi(u,v)= \varphi^{-1}( \min( \varphi(0), \varphi(u)+\varphi(v))),\quad (u,v)\in [0,1],
\end{align*}
for an Archimedean generator $\varphi:[0,1] \rightarrow [0,\infty)$ which is a convex, strictly decreasing and continuous function with $\varphi(1)=0$.
It is shown in~\citet{jaworski2004uniform} that,
if $\mathcal E_\varphi(0)=\lim_{x \downarrow 0}x \varphi'(x)/\varphi(x)=-\theta_0$, $0 < \theta_0 < \infty$,
then
\begin{align*}
\Lambda(u,v;C_\varphi)=uv  (u^{\theta_0} + v^{\theta_0})^{-1/\theta_0}.
\end{align*}
Therefore, we have that
\begin{align*}
\Lambda\left(b,\frac{1}{b};C_\varphi\right)=(b^{-\theta_0}+b^{\theta_0})^{-1/\theta_0}, \quad b \in (0,\infty).
\end{align*}
Since $\Lambda\left(b,1/b;C_\varphi\right)$ is maximized at $b^\ast=1$, we have that $\lambda(\Lambda)=\lambda^\ast(\Lambda)=2^{-1/\theta_0}$.
In addition, by~\eqref{eq:formula:maximal:atcm:max}, we have that $\overline \lambda(\Lambda)=\lim_{b\downarrow 0}(1+b^{2\theta_0})^{-1/\theta_0}=1$.
Examples of such Archimedean copulas include a Clayton copula with the generator $\varphi_\theta(t)=(t^{-\theta}-1)/\theta$, $\theta \in (0,\infty)$, which satisfies $\mathcal E_{\varphi_\theta}(0)=-\theta$; see Figure~\ref{fig:example:tail:dependence:functions} (i) for examples of the curve $b \mapsto \Lambda\left(b,1/b;C_{\varphi_\theta}\right)$.

It is also shown in~\citet{jaworski2004uniform} that, if  $\mathcal E_\varphi(1)=-\lim_{x \downarrow 0}x \varphi'(1-x)/\varphi(1-x)=\theta_1$, $1\le\theta_1<\infty$,
then the survival Archimedean copula has the tail copula
\begin{align*}
\Lambda(u,v; \hat C_\varphi)=u + v - (u^{\theta_1} + v^{\theta_1})^{1/\theta_1}.
\end{align*}
Therefore, we have that
\begin{align*}
\Lambda\left(b,\frac{1}{b}; \hat C_\varphi\right)=b + \frac{1}{b} - \left(b^{\theta_1} + b ^{-\theta_1}\right)^{1/\theta_1}.
\end{align*}
This function is also maximized at $b^\ast=1$, and thus we have that $\lambda(\Lambda)=\lambda^\ast(\Lambda)=2-2^{1/\theta_1}$.
In addition,~\eqref{eq:formula:maximal:atcm:max} yields $\overline \lambda(\Lambda)=1$.
Examples of such Archimedean copulas include a Gumbel copula with the generator $\varphi_\theta(t)=(-\log t)^{\theta}$, $\theta \in [1,\infty)$, which satisfies $\mathcal E_{\varphi_\theta}(1)=\theta$.
\end{example}

\begin{example}[Asymmetric Gumbel and Galambos copulas]\label{example:maximal:tail:dependence:ev:copulas}
The \emph{asymmetric Gumbel} and \emph{Galambos copulas} are EV copulas with the respective Pickands dependence functions given by
\begin{align*}
A_{\alpha,\beta,\theta}^\text{Gu}(w)&=(1-\alpha)w + (1-\beta)(1-w) + \{(\alpha w)^\theta + (\beta (1-w))^\theta\}^{1/\theta},\quad 1\leq \theta < \infty,\   0<\alpha,\beta\leq 1,\\
A_{\alpha,\beta,\theta}^\text{Ga}(w)&=1-\{(\alpha w)^{-\theta} + (\beta (1-w))^{-\theta}\}^{-1/\theta},\quad 0 <  \theta < \infty,\  0<\alpha,\beta\leq 1.
\end{align*}
Therefore, by~\eqref{eq:ev:copula:tail:dependence:function}, their survival copulas have the tail copulas
\begin{align*}
\Lambda_{\alpha,\beta,\theta}^\text{Gu}\left(b,\frac{1}{b}\right)&=
(\alpha - \beta) b + \beta \frac{b^2+1}{b}-b\left\{\alpha^\theta + \beta^\theta b^{-2\theta}\right\}^{1/\theta},\\
\Lambda_{\alpha,\beta,\theta}^\text{Ga}\left(b,\frac{1}{b}\right)&=
b \left\{\alpha^{-\theta} + \beta^{-\theta}b^{2\theta}\right\}^{-1/\theta},
\end{align*}
where $b \in (0,\infty)$; see Figure~\ref{fig:example:tail:dependence:functions}
(iii), (iv), (v), (vi) for examples of the corresponding curves.  A rather tedious
calculation shows that both tail copulas are maximized at $b^\ast=\sqrt{\beta/\alpha}$, and we obtain the formulas
\begin{align*}
\lambda(\Lambda_{\alpha,\beta,\theta}^{\text{Gu}})&=\alpha + \beta -(\alpha^\theta + \beta^\theta)^{1/\theta},\quad
\lambda(\Lambda_{\alpha,\beta,\theta}^{\text{Ga}})=(\alpha^{-\theta}+\beta^{-\theta})^{-1/\theta},\\
\lambda^\ast(\Lambda_{\alpha,\beta,\theta}^{\text{Gu}})&=\left(2-2^{1/\theta}\right)\sqrt{\alpha\beta},\quad
\lambda^\ast(\Lambda_{\alpha,\beta,\theta}^{\text{Ga}})=2^{-1/\theta}\sqrt{\alpha\beta},\\
\overline \lambda(\Lambda_{\alpha,\beta,\theta}^{\text{Gu}})&=\overline \lambda(\Lambda_{\alpha,\beta,\theta}^{\text{Ga}})=\max(\alpha,\beta).
\end{align*}
\end{example}

\section{Numerical experiments}\label{sec:numerical:experiments}

In this section we conduct numerical experiments to show the performance of the proposed tail concordance measures for various copulas.
The proposed measures are computed based on the empirical tail copula constructed from $n$ (pseudo-)observations of the underlying copula.
Namely, for a sequence $k=k(n)$ such that $k/n \rightarrow 0$ as $n \rightarrow \infty$, an empirical tail copula is constructed by $\tilde \Lambda^{[n,k]}(u,v)=\tilde C^{[n]}(ku/n,kv/n)/(k/n)$, $(u,v)\in \IR_{+}^2$, where $\tilde C^{[n]}$ is an empirical copula constructed from (pseudo-)observations.
Estimators of the proposed measures are then constructed by replacing $\Lambda$ with the empirical counterpart $\tilde \Lambda^{[n,k]}$.
The readers are referred to Appendix~\ref{sec:estimation:proposed:measures} for detailed construction and results on statistical inference of the proposed measures.

\subsection{Simulation study}\label{subsec:simulation:study}

\begin{figure}[t!]
  \centering
  \vspace{-25mm}
  \includegraphics[width=10 cm]{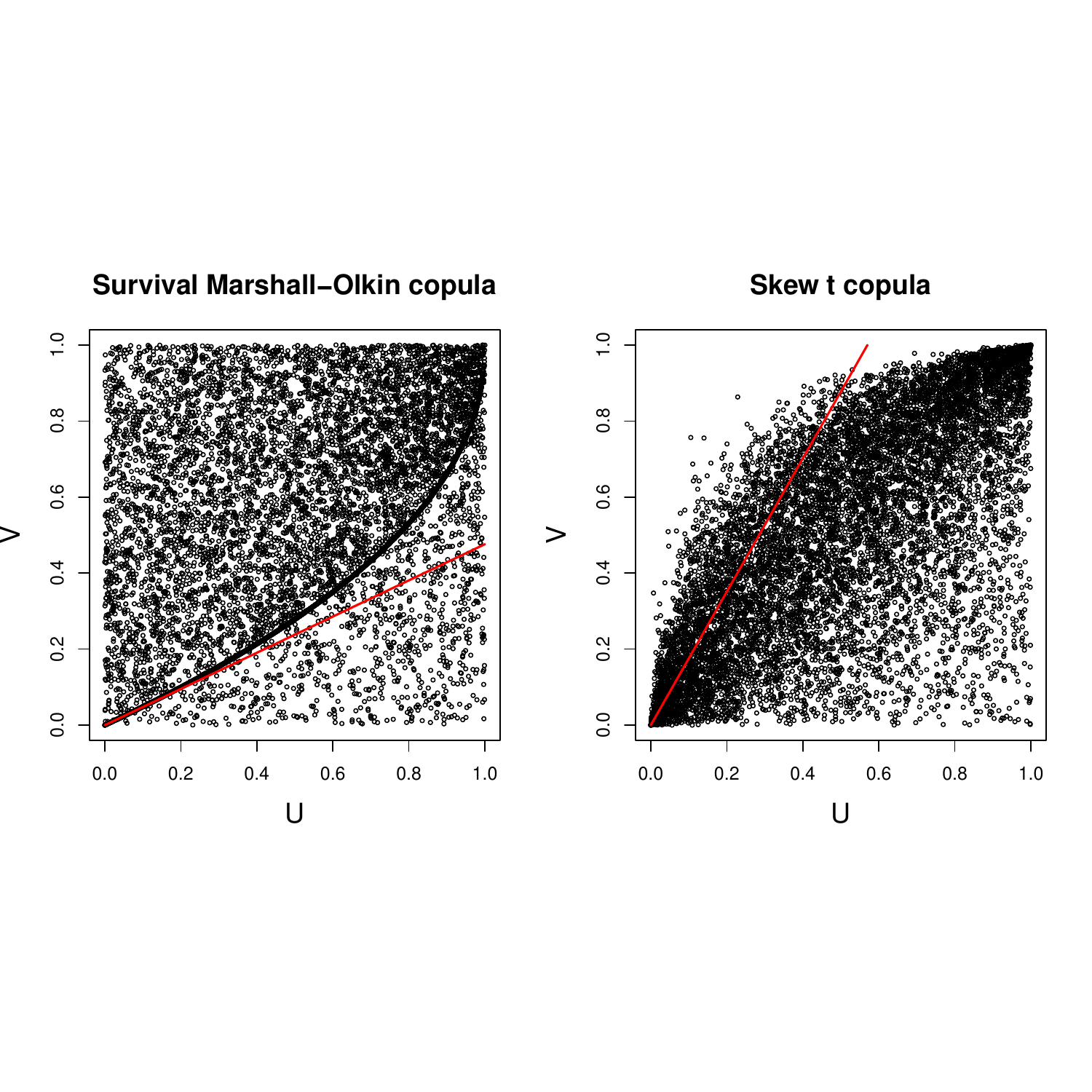}
  \vspace{-25mm}
  \caption{Scatter plots of the survival Marshall--Olkin copula
    $\hat C_{\alpha,\beta}^{\text{MO}}$ and the skew $t$ copula
    $C_{\nu,\delta_1,\delta_2,\gamma}^\text{ST}$ where the parameters are given
    by $(\alpha,\beta)=(0.353,0.75)$ and
    $(\nu,\delta_1,\delta_2,\gamma)=(5,0.8,-0.8,0.95)$. The sample size is
    $n=10^4$.  The red line indicates the estimated line $y=(1/b^\ast)^2 x$ such that $\Lambda(b^\ast,1/b^\ast)=\sup_{b \in (0,\infty)}\Lambda(b,1/b)$. }
  \label{fig:scatter:plot:SMO:skew:t:copulas}
\end{figure}

We first conduct a simulation study to estimate the proposed measures for two parametric copulas;
one is a survival Marshall--Olkin copula
$\hat C_{\alpha,\beta}^{\text{MO}}$ as in
Example~\ref{example:marshall:olkin:copula:maximal:tail:dependence:measure} for
$(\alpha,\beta)=(0.353,0.75)$, and another is a skew $t$ copula
$C_{\nu,\delta_1,\delta_2,\gamma}^\text{ST}$ for
$(\nu,\delta_1,\delta_2,\gamma)=(5,0.8,-0.8,0.95)$ as in
\cite{smith2012modelling}.
Figure~\ref{fig:scatter:plot:SMO:skew:t:copulas} shows scatter plots of these copulas.

\begin{figure}[t!]
  \centering
  \includegraphics[width=16 cm]{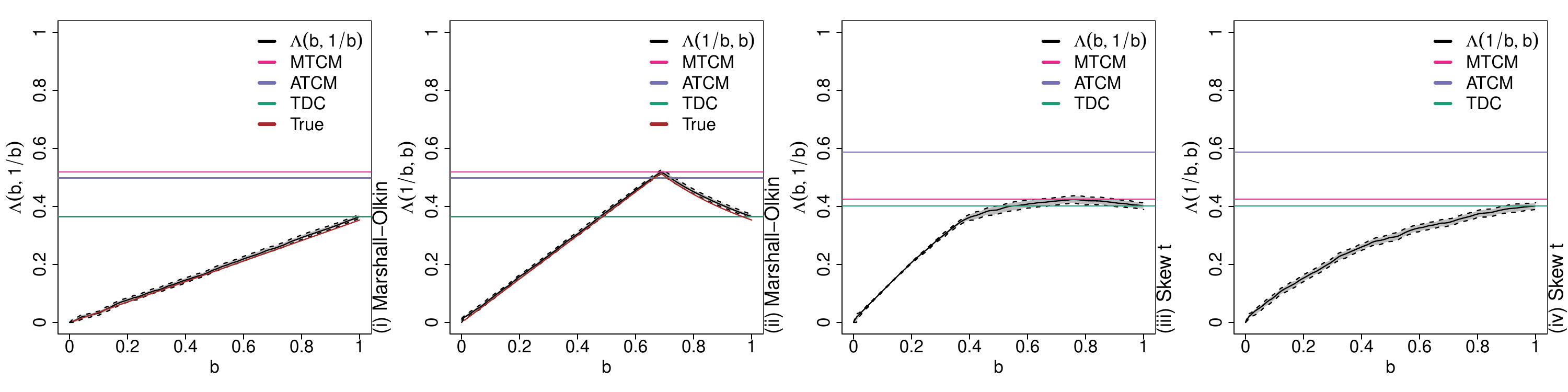}
  \caption{Plots of the estimates of the tail copulas
    $b \mapsto \Lambda\left(b,1/b\right)$ and
    $b \mapsto \Lambda\left(1/b,b\right)$ for (i), (ii) a survival
    Marshall--Olkin copula (left two plots); and (iii), (iv) a skew $t$ copula (right two plots) based on $n=10^6$
    samples.  The black solid lines indicate bootstrap means of the empirical
    tail copulas and the dashed lines are the 95\% bootstrap confidence intervals based on $B=100$ bootstrap replications.  The brown solid lines in (i) and (ii) are
    the true curves derived in
    Example~\ref{example:marshall:olkin:copula:maximal:tail:dependence:measure}.
    The colored horizontal lines represent the bootstrap estimates of the TCMs.
}
  \label{fig:Lambda:tail copulas:SMO:skew:t:copulas}
\end{figure}

We first simulate $n=10^6$ samples from $\hat C_{\alpha,\beta}^{\text{MO}}$ and $C_{\nu,\delta_1,\delta_2,\gamma}^\text{ST}$, and then evaluate the corresponding empirical tail copulas at $(b,1/b)$ and $(1/b,b)$ for $b \in \left\{1/L,\dots,1\right\}$, where $L=100$.
The threshold $k \in \IN$ is chosen to be $k_\text{SMO}=0.015\,\times\,N = 15000$ for the survival Marshall--Olkin copula and $k_\text{ST} =  0.005 \,\times\,N = 5000$ for the skew $t$ copula.
These values are determined by a graphical plateau--finding search where the estimates of the TDC are plotted against various values of $k$, and an optimal $k$ is chosen from an interval on which the estimated TDCs are roughly constant; see \cite{schmidt2006non} for details.
Finally, based on the empirical tail copulas, we estimate the TDC $\lambda$, the discrete uniform ATCM $\lambda_{\text{U}_L}$, where the equal weights are put on $(b,1/b)$ and $(1/b,b)$ for $b \in \left\{1/L,\dots,1\right\}$, and the MTCM $\lambda^\ast$ with the estimators provided in Section~\ref{sec:estimation:proposed:measures}.
To reduce the computational cost, we maximize the empirical tail copulas only over $\{1/L,\dots,1,L/(L-1),\dots,L\}$ when estimating the MTCM.
Note that the true tail copulas and TCMs are available for survival Marshall--Olkin copulas; see Example~\ref{example:marshall:olkin:copula:maximal:tail:dependence:measure}.
As mentioned in~\cite{smith2012modelling}, the analytical calculation of the TDC for skew $t$ copulas may not be straightforward.
In Figure~\ref{fig:Lambda:tail copulas:SMO:skew:t:copulas} and Table~\ref{Table:estimates:tail:dependence:measures}, we report the estimates of the tail copulas and the TCMs with the 95\% bootstrap confidence intervals (CIs) based on $B=100$ bootstrap replications.
As in Remark~\ref{remark:attainability:mtcm}, we also report the maximizer $b^\ast$ of the MTCM as a measure of tail non-exchangeability.
To treat $b^\ast$ and $1/b^\ast$ on the same scale, we report
\begin{align}\label{eq:def:b:angle:bracket}
\langle b^\ast \rangle =\begin{cases}
b^\ast, & \text{ if } 0<b^\ast \leq 1,\\
2-1/b^\ast \in [1,2), & \text{ if } 1<b^\ast,\\
\end{cases}
\end{align}
instead of $b^\ast$.
As a result, the square attains the MTCM if $\langle b^\ast \rangle=1$, and a more elongated rectangle attains the MTCM as $\langle b^\ast \rangle$ goes to $0$ or $2$.

\begin{table}[t!]
  \centering
\begin{small}
    \begin{tabular}{
    l
    cccccc
    }
    \toprule
       &  & $\lambda$ &  $\lambda^\ast$ & $\lambda_{\text{U}_L}$ & $\langle b^\ast\rangle $ \\
    \midrule
   \multicolumn{6}{l}{(1) Survival Marshall--Olkin copula $\hat C_{\alpha,\beta}^{\text{MO}}$}\\[4pt]

&Estimate & 0.367 &  {\bf 0.523} &0.504 & 1.310 \\
&  95\% CI & (0.359,\,0.374) &  (0.514,\,0.529) & (0.493,\,0.512) &(1.300,\,1.320)  \\

   \multicolumn{6}{l}{(2) Skew $t$ copula $C_{\nu,\delta_1,\delta_2,\gamma}^\text{ST}$}\\[4pt]

&  Estimate& 0.407 & 0.417 &  {\bf 0.581} &0.779 \\
&    95\% CI & (0.396,\,0.416) & (0.407,\,0.430) & (0.565,\,0.596) & (0.650,\,0.895) \\
      \bottomrule
    \end{tabular}
    \end{small}
    \caption{Bootstrap estimates and 95\% CIs of the TDC $(\lambda)$, the discrete uniform ATCM ($\lambda_{\text{U}_L}$), the MTCM ($\lambda^\ast$) and its maximizer $\langle b^\ast \rangle$ in~\eqref{eq:def:b:angle:bracket} of the survival Marshall--Olkin copula and the skew $t$ copula with $B=100$ bootstrap replications.}
    \label{Table:estimates:tail:dependence:measures}
\end{table}

From Figure~\ref{fig:Lambda:tail copulas:SMO:skew:t:copulas} and Table~\ref{Table:estimates:tail:dependence:measures}, we observe that the bootstrap CIs for the tail copulas and the proposed TCMs are sufficiently narrow.
In particular, the estimated curves of $b \rightarrow \Lambda(b,1/b)$ are sufficiently close to the true ones in Figures~\ref{fig:Lambda:tail copulas:SMO:skew:t:copulas} (i) and (ii).
As we can see, both the survival Marshall--Olkin and the skew $t$ copula have non-exchangeable tail dependence.
The two copulas, however, have different features of tail dependence since the tail copula in Figure~\ref{fig:Lambda:tail copulas:SMO:skew:t:copulas} (ii) has a sharp kink around $b=0.7$ whereas the tail copula in Figure~\ref{fig:Lambda:tail copulas:SMO:skew:t:copulas} (iii) constantly takes on large values for $0.6 \leq b \leq 1$.
Since the former feature is captured by the MTCM and the latter one by the ATCM, the gap between the TDC and the MTCM is large for the survival Marshall--Olkin copula and the ATCM of the skew $t$ copula is larger than that of the survival Marshall--Olkin copula.

\subsection{Real data analysis}\label{subsec:real:data:analysis}

\begin{figure}[t!]
  \centering
  \includegraphics[width=16 cm]{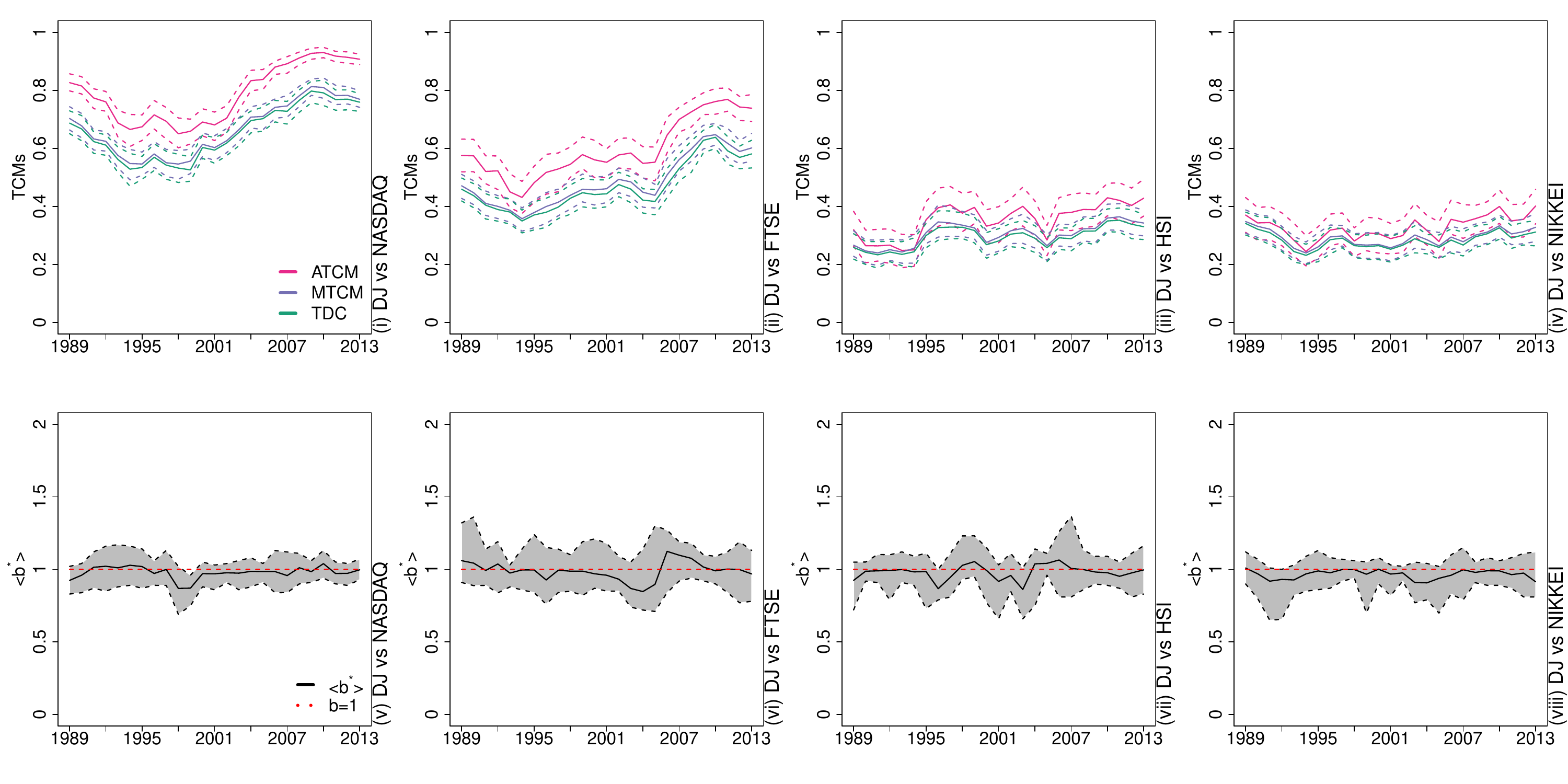}
  \caption{
  Plots of the estimates of the proposed TCMs for filtered stock returns of (i), (v)
    (DJ, NASDAQ); (ii), (vi) (DJ, FTSE); (iii), (vii) (DJ, HSI); and (iv), (viii) (DJ,
    NIKKEI).
    The solid lines indicate bootstrap means and the dashed lines are the 90\% bootstrap confidence intervals based on $B=100$ bootstrap replications.}
  \label{fig:LambdaF:TS}
\end{figure}

To investigate financial interconnectedness between different countries in times of a stressed economy, we compare the proposed TCMs for the returns of the stock indices DJ, NASDAQ, FTSE, HSI and NIKKEI from 1987-01-05 to 2015-12-30.
We particularly focus on the relationship of the DJ to the other indices to compare the corresponding domestic and international relationships within and to the US market, respectively.
To this end, a rolling window analysis is conducted on which, for each year from 1989 to 2013, the last and next two years of data are included in the window, resulting in 5 year windows with sample size $n\approx 250 \,\times\, 5 = 1250$ each.
For each window, we filter the marginal return series by GARCH(1,1) models with skew $t$ innovations.
The residuals are then rank-transformed to obtain the pseudo--observations of the underlying copula.
Based on the pseudo--observations, we conduct the same analysis as in Section~\ref{subsec:simulation:study}, where $L=100$ and the threshold $k$ is chosen to be $0.2$ times the sample size following~\citet[Section~5]{bormann2020detecting}.
The results are summarized in Figure~\ref{fig:LambdaF:TS}.

In Figure~\ref{fig:LambdaF:TS} (i)--(iv), we observe that the overall level of tail dependence is high for the pairs (DJ,\,NASDAQ) and (DJ,\,FTSE), and is relatively low for the pairs (DJ,\,HSI) and (DJ,\,NIKKEI).
For all the pairs, the gap between the TDC and the MTCM is not significant.
For the cases (i) and (ii), an increasing trend of tail dependence can be observed since around the financial crisis of 2007--2009.
Such a trend is barely observable for the cases (iii) and (iv).
In Figure~\ref{fig:LambdaF:TS} (v)--(viii), the maximizer $\langle b^\ast\rangle$ of the MTCM fluctuates over time around $1$ for all the cases, but some bumps are visible during financial crises, such as
the collapse of the Japanese asset price bubble in 1992 as seen in Figure~\ref{fig:LambdaF:TS} (viii),
the Asian and millennium crisis accumulating into the Dot-Com crisis in 1995--2003 as seen in Figure~\ref{fig:LambdaF:TS} (v), (vii) and (viii), and
the subprime mortgage crisis in 2007--2009 as seen in Figure~\ref{fig:LambdaF:TS} (vi), (vii) and (viii).
The comparably small fluctuation in Figure~\ref{fig:LambdaF:TS} (v) implies tail exchangeability between two US stock indices.

\section{Conclusion}\label{sec:conclusion}

We proposed two novel tail dependence measures which we call the \emph{maximal tail concordance measure (MTCM)} and the $\mu$-\emph{average tail concordance measure (ATCM)}.
Both measures are based on tail copulas and possess clear probabilistic interpretations.
The MTCM evaluates the largest limiting tail probability over all comparable rectangles in the tail, and the ATCM is a normalized average of these limiting tail probabilities.
With these interpretations, the MTCM is useful for extracting the most distinctive feature of tail dependence, and ATCMs can be applied when specific joint tail events are of particular importance from a practical point of view.
We showed that the two measures satisfy axiomatic properties naturally required for quantifying tail dependence.
In contrast to the TDC, both proposed measures also capture non-exchangeable tail dependence.
The choice of the angular measure $\mu$ of the ATCM was addressed via practical examples.
Bounds of $\mu$-ATCMs were also explored.
We showed that the minimal ATCM is the TDC, and that, unlike the MTCM, the maximal ATCM does not extract informative features of tail dependence.
Furthermore, we provided analytical formulas for various parametric copulas and simulation studies to support the use of the proposed measures.
Finally, a real data analysis revealed tail dependence and tail non-exchangeability of the return series of stock indices, particularly in periods of financial distress.

Further investigation is needed for statistical inference of the proposed measures and for the relationship between the MTCM and the tail indices proposed by \citet{furman2015paths} and \citet{genest2021class}.
Moreover, multivariate extensions and compatibility problems of the proposed TCMs are also interesting avenues of future research; see \cite{frahm2006extremal}, \cite{schmid2007multivariateTail}, \cite{li2009orthant} and \cite{gijbels2020multivariate} for the former problem and \cite{embrechtshofertwang2016} and \cite{hofert2019compatibility} for the latter.

\section*{Acknowledgements}
We would like to thank insightful comments from the editor and two anonymous reviewers.
We are also grateful to Ricardas Zitikis, Ruodu Wang and Shunichi Nomura for their valuable comments.
Takaaki Koike was supported by Japan Society for the Promotion of Science (JSPS KAKENHI Grant Number JP21K13275).
Shogo Kato acknowledges financial support from JSPS KAKENHI Grant Number JP20K03759.
Marius Hofert would like to thank the financial support from the Natural Sciences and Engineering Research Council of Canada (NSERC RGPIN-2020-04897 and RGPAS-2020-00093).

\section*{Competing interests}
The authors declare none.

\bibliographystyle{apalike}

\appendix

\section*{Appendices}

\section{Basic properties of tail copulas}\label{sec:basic:properties:tail:copulas}

Tail copulas play an important role for quantifying extremal co-movements between random variables.
In this section we summarize basic properties of tail copulas.
To this end, let $\mathcal L$ be the set of all tail copulas.

\begin{proposition}[Basic properties of tail copulas]
\label{prop:basic:properties:tail:dependence:function}
Let $\Lambda\in \mathcal L$.
\begin{enumerate}[label=\arabic*), itemsep=0pt]
    \item\label{basic:prop:item:2:increasing} {\bf ($2$-increasingness)}: $\Lambda$ is 2-increasing, that is, $\Lambda(u',v')-\Lambda(u',v)-\Lambda(u,v')+\Lambda(u,v)\geq 0$ for every $0\leq u\leq u' $ and $0\leq v\leq v'$.
    \item\label{basic:prop:item:monotonicity} {\bf (Monotonicity)}: $\Lambda(u,v)\leq \Lambda(u',v')$  for every $0\leq u\leq u'$ and $0\leq v\leq v'$. If $\Lambda\not \equiv 0$, then $\Lambda(u,v)<\Lambda(u',v')$ for every $0<u < u'$ and $0<v < v'$.
    \item\label{basic:prop:item:grounded} {\bf (Groundedness)}: $\Lambda$ is grounded, that is, $\Lambda(u,v)=0$ if $u=0$ or $v=0$.
    \item\label{basic:prop:item:one:homogeneity} {\bf (Positive homogeneity)}: $\Lambda(tu,tv)=t\Lambda(u,v)$ for every $t \geq 0$ and $(u,v)\in\IR_{+}^2$.
    \item\label{basic:prop:item:degeneracy} {\bf (Degeneracy)}: $\Lambda(u,v)= 0$ for all $(u,v)\in \IR_{+}^2$ if and only if $\Lambda(u_0,v_0)=0$ for some $u_0,v_0>0$.
    \item\label{basic:prop:item:coherence} {\bf (Coherence)}: If $C\preceq C'$ for $C,C' \in \mathcal C_2^\text{L}$, then  $\Lambda(u,v;C)\leq \Lambda(u,v;C')$ for all $(u,v)\in \IR_{+}^2$.
    \item\label{basic:prop:item:bounds} {\bf (Bounds)}:
    $0\leq \Lambda(u,v)\leq \min(u,v)$ for all $(u,v)\in\IR_{+}^2$ and the bounds are attainable.
    \item\label{basic:prop:item:homo:bounds} {\bf (Max-min inequalities)}:  $(s\wedge t)\Lambda(u,v)\leq \Lambda(su,tv)\leq (s\vee t)\Lambda(u,v)$ for every $s,t\geq 0$ and $(u,v)\in \IR_{+}^2$.
\end{enumerate}
\end{proposition}

\begin{proof}
\ref{basic:prop:item:2:increasing},~\ref{basic:prop:item:grounded},~\ref{basic:prop:item:one:homogeneity} and~\ref{basic:prop:item:degeneracy} can be found in \citet[Propositions~2.1 and 2.2]{joe2010tail}.
\ref{basic:prop:item:monotonicity} is shown in \citet[Theorems~1 and 2]{schmidt2006non}.
\ref{basic:prop:item:coherence} follows directly from the definition of $\Lambda$.
\ref{basic:prop:item:bounds} is implied by~\ref{basic:prop:item:coherence}.
\ref{basic:prop:item:homo:bounds} follows from~\ref{basic:prop:item:monotonicity} and~\ref{basic:prop:item:one:homogeneity}.
\end{proof}

\begin{proposition}[Continuity and derivatives of tail copulas]
\label{prop:continuity:derivatives}
Let $\Lambda\in \mathcal L$.
\begin{enumerate}[label=\arabic*), itemsep=0pt]
    \item\label{basic:prop:item:continuity} {\bf (Continuity)}: $|\Lambda(u,v)-\Lambda(u',v')|\leq |u-u'|+|v-v'|$ for every $(u,v),(u',v') \in \IR_{+}^2$, and thus $\Lambda$ is Lipschitz continuous.
    \item\label{basic:prop:item:partial:derivatives} {\bf (Partial derivatives)}:  The partial derivatives $\mathrm{D}_1\Lambda(u,v)$ and $\mathrm{D}_2\Lambda(u,v)$ exist almost everywhere on $(u,v)\in\IR_{+}^2$. Moreover, $0\leq \mathrm{D}_1\Lambda(u,v),\ \mathrm{D}_2\Lambda(u,v)\leq 1$ and the functions
    \begin{align*}
    v\mapsto \mathrm{D}_1\Lambda(u,v)\quad \text{and}\quad
    u\mapsto \mathrm{D}_2\Lambda(u,v)
    \end{align*}
    are increasing almost everywhere on $\IR_{+}$.
    \item\label{basic:prop:item:euler:formula}  {\bf (Euler's theorem)}: $\Lambda(u,v)=u\mathrm{D}_1\Lambda(u,v)+v \mathrm{D}_2 \Lambda(u,v)$ for every $(u,v)\in\IR_{+}^2$.
\end{enumerate}
\end{proposition}

\begin{proof}
\ref{basic:prop:item:continuity} and \ref{basic:prop:item:partial:derivatives} are shown in~\citet[Theorems~1 and~3]{schmidt2006non}, respectively.
\ref{basic:prop:item:euler:formula} is the well-known Euler's homogeneous function theorem.
\end{proof}

\section{Proofs}\label{sec:proofs}

\subsection{Proof of Proposition~\ref{prop:tail:concordance:order:tail:probabilities}}

\noindent
\ref{item:tail:concordance:order:theta:angle}
Let us write a point $(u,v)\in (0,\infty)^2$ in polar coordinates via $(u,v)=r(\cos \theta, \sin \theta)$ for $r > 0$ and $\theta \in (0,\pi/2)$.
By positive homogeneity of $\Lambda$, we have that
\begin{align*}
\Lambda(u,v)=\frac{r}{r_\theta}\Lambda\left(r_\theta \cos(\theta),r_\theta \sin(\theta)\right),
\end{align*}
and thus $\Lambda \preceq \Lambda'$ is equivalent to~\eqref{eq:simplified:condition:Lambda:concordance:order:theta:angle}.

\noindent
\ref{item:tail:concordance:order:tail:probability}
Notice that $\Lambda(u,1)=\lim_{p\downarrow 0}\Prob(U\leq pu\,|\, V \leq p)$ and $\Lambda(1,v)=\lim_{p\downarrow 0}\Prob(V\leq pv\,|\, U \leq p)$ for $u,v\in (0,1]$.
Therefore, the desired equivalence holds from Part~\ref{item:tail:concordance:order:theta:angle} of this proposition by taking $r_\theta=1/\max(\cos\theta,\sin\theta)$.

\subsection{Proof of Proposition~\ref{prop:representation:mu:atcm}}

Representation~\eqref{eq:mu:atcm:representation} is obtained by taking $w_0=\mu(\{1\})$, $w_1=\mu((0,1))$, $w_2 = \mu((1,\infty))$ and $\mu_1(A) = \mu (A\cap (0,1))/w_1$ and $\mu_2(A) = \mu (\{1/b;b \in A\}\cap (0,1))/w_2$, $A \in \mathfrak B ((0,\infty))$, if $w_1>0$ and $w_2>0$, respectively.

\subsection{Proof of Proposition~\ref{prop:axiomatic:properties:proposed:tcms}}

\noindent
\ref{prop:item:mtcm:is:convex:tcm}
The first three axioms and convexity are straightforward to check, and thus
it remains to show the continuity axiom and the strictness.
Suppose that $\Lambda^{[n]} \in \mathcal L$, $n=1,2,\dots,$ converge to $\Lambda \in \mathcal L$ pointwise as $n\rightarrow \infty$.
Then the convergence is uniform by Proposition~\ref{prop:continuity:derivatives}~\ref{basic:prop:item:continuity},
Therefore, we have that
\begin{align*}
\lim_{n\rightarrow \infty}\lambda^\ast(\Lambda^{[n]})=\lim_{n\rightarrow \infty}\max_{b \in (0,\infty)}\Lambda^{[n]}\left(b,\frac{1}{b}\right)=\max_{b \in (0,\infty)}\Lambda\left(b,\frac{1}{b}\right)=\lambda^\ast(\Lambda),
\end{align*}
where the supremum in $\lambda^\ast$ can be replaced by the maximum as mentioned in Remark~\ref{remark:attainability:mtcm}.
To show the strictness of $\lambda^\ast$, we first show the following corollary.
\begin{corollary}[Tail independence and tail comonotonicity for the TDC]
\label{cor:tail:independence:tail:comonotonicity:tdc}
A tail copula $\Lambda \in \mathcal L$ is tail independent if and only if $\lambda(\Lambda)=0$, and $\Lambda$ is tail comonotonic if and only if $\lambda(\Lambda)=1$.
\end{corollary}

\begin{proof}
The equivalence in the tail independent case follows directly from Proposition~\ref{prop:basic:properties:tail:dependence:function}~\ref{basic:prop:item:degeneracy}.
If $\Lambda$ is tail comonotonic, then $\lambda(\Lambda)=\lambda(\overline \Lambda)=1$.
To show the converse, suppose by way of contradiction that $\lambda(\Lambda)=1$ but that there exists $(u_0,v_0)\in \IR_{+}^2$ such that $\Lambda(u_0,v_0)\neq \overline \Lambda(u_0,v_0)$.
Note that $u_0,v_0>0$, otherwise $\Lambda(u_0,v_0)= \overline \Lambda(u_0,v_0)=0$.
By~\eqref{eq:inequality:tail:dependence:coefficient}, we have that
\begin{align*}
1=\lambda(\Lambda)\leq \lambda_{\delta_{\sqrt{u_0/v_0}}}(\Lambda)=\frac{\Lambda(u_0,v_0)}{\overline \Lambda(u_0,v_0)}<1,
\end{align*}
which is a contradiction.
\end{proof}

By Corollary~\ref{cor:tail:independence:tail:comonotonicity:tdc}, it suffices to show that $\lambda^\ast(\Lambda)=1$ implies $\lambda(\Lambda)=1$.
Suppose $\lambda^\ast(\Lambda)=1$.
As mentioned in Remark~\ref{remark:attainability:mtcm}, there exists $b^\ast \in (0,\infty)$ such that $\lambda^\ast(\Lambda)=\Lambda\left(b^\ast,1/b^\ast\right)$.
Suppose that $0<b_\ast<1$.
Then, by Proposition~\ref{prop:basic:properties:tail:dependence:function}~\ref{basic:prop:item:bounds}, we have that
\begin{align*}
1=\Lambda\left(b^\ast,\frac{1}{b^\ast}\right)\leq \overline \Lambda\left(b^\ast,\frac{1}{b^\ast}\right)=b_\ast
\end{align*}
and thus that $1\leq b_\ast$, which is a contradiction.
Similarly, if we assume that $1<b_\ast$, then we have that $1=\Lambda\left(b_\ast,1/b_\ast\right)\leq M\left(b_\ast,1/b_\ast\right)=1/b_\ast$ and thus that $b_\ast\leq 1$, which is again a contradiction.
Therefore, we have that $b_\ast=1$, which yields $\lambda(\Lambda)=\Lambda(1,1)=\Lambda\left(b_\ast,1/b_\ast\right)=\lambda^\ast(\Lambda)=1$.

\noindent
\ref{prop:item:atcm:is:linear:tcm}
Linearity and Axioms~\ref{def:axioms:item:monotonicity} and~\ref{def:axioms:item:continuity} immediately follow by definition of $\lambda_\mu$.
Axioms~\ref{def:axioms:item:normalization} and~\ref{def:axioms:item:tail:independence} are straightforward to show by~\eqref{eq:inequality:tail:dependence:coefficient}
and Corollary~\ref{cor:tail:independence:tail:comonotonicity:tdc}.

\noindent
\ref{prop:item:condition:atcm:is:strict:tcm}
By Corollary~\ref{cor:tail:independence:tail:comonotonicity:tdc}, it suffices to show that $\lambda_\mu(\Lambda)=1$ implies $\lambda(\Lambda)=1$.
Suppose, by way of contradiction, that $\lambda=\lambda(\Lambda)\leq 1- \delta<1$ for some $\delta>0$.
By Condition~\eqref{eq:condition:probability:around:diagonal}, at least one of the following two cases holds: $\mu((1-\epsilon,1])>0$ for any $\epsilon>0$ or $\mu([1,1+\epsilon))>0$ for any $\epsilon>0$.
Assume that the former case is fulfilled; the latter case can be shown similarly.
Let
\begin{align*}
c(\delta)=\frac{1-\delta}{1-\delta/2}\in (0,1)\quad\text{and}\quad
R_{\delta}=\left[\sqrt{c(\delta)},1\right].
\end{align*}
Then $1/b \leq b/c(\delta)$ for $b \in R_\delta$.
Together with $\Lambda(1,1)=\lambda(\Lambda)\leq 1-\delta$, we have that, for $b \in R_\delta$,
\begin{align*}
\Lambda\left(b,\frac{1}{b}\right)&\leq \Lambda\left(\frac{1}{b},\frac{1}{b}\right)=\frac{1}{b}\Lambda(1,1)\leq \frac{1}{c(\delta)}b(1-\delta)=\left(1-\frac{\delta}{2}\right)b=\left(1-\frac{\delta}{2}\right)\overline \Lambda\left(b,\frac{1}{b}\right).
\end{align*}
Since $\Lambda(b,1/b)\leq \overline \Lambda(b,1/b)$ for $b \in (0,\infty)\backslash R_\delta$, we have that
\begin{align*}
\lambda_\mu(\Lambda)&=\frac{\int_{(0,\infty)}\Lambda(b,1/b) \rd \mu(b)}{\int_{(0,\infty)} \overline \Lambda(b,1/b) \rd \mu(b)}\leq \frac{\int_{(0,\infty)} \overline \Lambda(b,1/b) \rd \mu(b)-\frac{\delta}{2}\int_{R_\delta}\overline \Lambda(b,1/b)\rd \mu(b)}{\int_{(0,\infty)} \overline \Lambda(b,1/b) \rd \mu(b)}\\
&=1-\frac{\delta}{2}\frac{\int_{R_\delta}\overline \Lambda(b,1/b)\rd \mu(b)}{\int_{(0,\infty)} \overline \Lambda(b,1/b)\rd \mu(b)}.
\end{align*}
Moreover, we have that
\begin{align*}
\int_{R_\delta}\overline \Lambda(b,1/b)\rd \mu(b)=\int_{R_\delta}b\rd \mu(b)\geq \sqrt{c(\delta)}\mu(R_\delta)>0,
\end{align*}
and hence $\lambda_\mu(\Lambda)<1$, which contradicts the assumption that $\lambda_\mu(\Lambda)=1$.

\subsection{Proof of Theorem~\ref{thm:characterization:minimal:maximal:atcm}}

\noindent
\ref{thm:item:characterization:minimal:TCM}
Inequalities \eqref{eq:inequality:tail:dependence:coefficient} directly follow from Proposition~\ref{prop:basic:properties:tail:dependence:function}~\ref{basic:prop:item:homo:bounds}.

\noindent
\ref{thm:item:characterization:maximal:average:tail:concordance:measure}
To show Equation~\eqref{eq:formula:maximal:atcm}, the inequality $\overline \lambda(\Lambda)\leq \sup_{b \in (0,\infty)}
\Lambda\left(b,1/b\right)
/
\overline \Lambda\left(b,1/b\right)
$ follows since
\begin{align*}
\int_{(0,\infty)}\Lambda(b,1/b)\rd \mu(b)=\int_{(0,\infty)}\overline \Lambda(b,1/b)\,\frac{\Lambda(b,1/b)}{\overline \Lambda(b,1/b)}\rd \mu(b)
\leq \sup_{b \in (0,\infty)}\frac{
\Lambda\left(b,1/b\right)
}{
\overline \Lambda\left(b,1/b\right)
}\int_{(0,\infty)}\overline \Lambda(b,1/b)\rd \mu(b).
\end{align*}
Now let $\Lambda^\star(b)=\Lambda(b,1/b)/\overline \Lambda(b,1/b)$.
By Proposition~\ref{prop:basic:properties:tail:dependence:function}~\ref{basic:prop:item:monotonicity} and Proposition~\ref{prop:continuity:derivatives}~\ref{basic:prop:item:continuity}, the map $b \mapsto \Lambda^\star(b)$ is continuous, decreasing on $(0,1]$ and increasing on $[1,\infty)$ almost everywhere.
Therefore, for any $\epsilon>0$, there exists $b^\star \in (0,\infty)$ such that $\sup_{(0,\infty)}\Lambda^\star(b)-\epsilon < \Lambda^\star(b^\star)$.
Since $\Lambda^\star(b^\star)=\lambda_{\delta_{b^\star}}(\Lambda)$ for $\delta_{b^\star} \in \mathcal M$, we obtain Equation~\eqref{eq:formula:maximal:atcm}.
Equality~\eqref{eq:formula:maximal:atcm:max} immediately follows from the monotonicity property of $b \mapsto \Lambda^\star(b)$.
Since $\delta_{0}$ and $\delta_{\infty}$ are not in $\mathcal M$, the supremum in~\eqref{eq:formula:maximal:atcm} is not attainable in general.

\noindent
\ref{thm:item:bounds:mtcm}
The inequality $\lambda(\Lambda)\leq\lambda^\ast(\Lambda)$ follows since
$\lambda(\Lambda)=\Lambda(1,1)\leq \sup_{b \in (0,\infty)}\Lambda\left(b,1/b\right)=\lambda^\ast(\Lambda)$.
As mentioned in Remark~\ref{remark:attainability:mtcm}, there exists $b^\ast \in (0,\infty)$ such that $\lambda^\ast(\Lambda)=\Lambda(b^\ast,1/b^\ast)$.
Since $\overline \Lambda(b^\ast,1/b^\ast)\leq 1$, we have, by~\eqref{eq:formula:maximal:atcm}, that
\begin{align*}
\lambda^\ast(\Lambda)=\Lambda\left(b^\ast,\frac{1}{b^\ast}\right)\leq\frac{\Lambda(b^\ast,1/b^\ast)}{
\overline \Lambda(b^\ast,1/b^\ast)}\leq \overline \lambda(\Lambda),
\end{align*}
which completes the proof.

\section{Estimation of the proposed measures}\label{sec:estimation:proposed:measures}

In this section we construct non-parametric estimatiors of the proposed TCMs based on the asymptotic results of empirical tail copulas developed in \cite{schmidt2006non} and \cite{bucher2013multiplier}.
To this end, let $(X_i,Y_i)$, $i=1,\dots,n$, $n \in \IN$, be an i.i.d.\ sample from a joint distribution function $H$ with continuous marginal distributions $X\sim F$ and $Y \sim G$, and copula $C \in \mathcal C_2^\text{L}$.
If $F$ and $G$ are known, then $(U_i,V_i)=(F(X_i),G(Y_i))$, $i=1,\dots,n$, is an i.i.d.\ sample from $C$.
If $F$ and $G$ are unknown, then
\begin{align*}
(\hat U_i,\hat V_i)=(\hat F^{[n]}(X_i),\hat G^{[n]}(Y_i))\quad \text{where}\quad \hat F^{[n]}(x)=\frac{1}{n+1}\sum_{i=1}^n \bone_{\{X_i \leq x\}},\
 \hat G^{[n]}(x)=\frac{1}{n+1}\sum_{i=1}^n \bone_{\{Y_i \leq x\}},
\end{align*}
are the \emph{pseudo--observations} from $C$.
Denote by
\begin{align*}
\tilde C^{[n]}(u,v)=\frac{1}{n}\sum_{i=1}^n \bone_{\{U_i\leq u,V_i\leq v\}}\quad \text{and}\quad
\hat C^{[n]}(u,v)=\frac{1}{n}\sum_{i=1}^n \bone_{\{\hat U_i\leq u,\hat V_i\leq v\}}
\end{align*}
the \emph{empirical copulas} based on the samples $(U_i,V_i)$ and $(\hat U_i,\hat V_i)$, respectively.
Define the \emph{empirical tail copulas}
\begin{align*}
\tilde \Lambda^{[n,k]}(u,v)=\frac{n}{k}\tilde C^{[n]}\left(\frac{ku}{n},\frac{kv}{n}\right)\quad\text{and}\quad
\hat \Lambda^{[n,k]}(u,v)=\frac{n}{k}\hat C^{[n]}\left(\frac{ku}{n},\frac{kv}{n}\right),
\end{align*}
where $k=k(n)$ is such that $k(n) \rightarrow \infty$ and $k=\operatorname{o}(n)$ as $n \rightarrow \infty$.

For $n,k \in \IN$ and $L>1$, define the estimators of $\lambda^\ast$ by
\begin{align}\label{eq:estimators:mtcm}
{\tilde \lambda}^{\ast[n,k,L]}=\max_{b\in[1/L,L]}\tilde \Lambda^{[n,k]}\left(b,\frac{1}{b}\right)\quad\text{and}\quad
{\hat \lambda}^{\ast[n,k,L]}=\max_{b\in[1/L,L]}\hat \Lambda^{[n,k]}\left(b,\frac{1}{b}\right)
\end{align}
and the estimators of $\lambda_\mu$, $\mu \in \mathcal M$, by
\begin{align}
\label{eq:estimators:mu:atcm:known:margin}
\tilde \lambda_{\mu}^{[n,k]}
 &=\frac{\int_{(0,\infty)}\tilde \Lambda^{[n,k]}(b,1/b)\rd \mu(b)}{\int_{(0,\infty)}\overline \Lambda(b,1/b)\rd \mu(b)}
 =\frac{1}{\int_{(0,\infty)}\overline \Lambda(b,1/b)\rd \mu(b)}
\,\frac{1}{k}
\sum_{i=1}^n \mu\left(
\left[
\frac{nU_i}{k},
\frac{k}{n V_i}
\right]
\right),\\
\label{eq:estimators:mu:atcm:unknown:margin}
\hat \lambda_{\mu}^{[n,k]}
 &=\frac{\int_{(0,\infty)}\hat \Lambda^{[n,k]}(b,1/b)\rd \mu(b)}{\int_{(0,\infty)}\overline \Lambda(b,1/b)\rd \mu(b)}
=\frac{1}{\int_{(0,\infty)}\overline \Lambda(b,1/b)\rd \mu(b)}
\,\frac{1}{k}
\sum_{i=1}^n \mu\left(
\left[
\frac{n{\hat U}_i}{k},
\frac{k}{n {\hat V}_i}
\right]
\right).
\end{align}

Since $\tilde \Lambda^{[n,k]}$ and $\hat \Lambda^{[n,k]}$ are step functions, the estimators in~\eqref{eq:estimators:mtcm} can be computed by maximizing $\tilde \Lambda^{[n,k]}$ and $\hat \Lambda^{[n,k]}$ over a finite number of points in $[1/L,L]^2$.
The estimators in~\eqref{eq:estimators:mu:atcm:known:margin} and~\eqref{eq:estimators:mu:atcm:unknown:margin} can be computed under the assumptions that
$\mu([s,t])$, $0<s\le t <\infty$, and the denominator $\int_{(0,\infty)}\overline \Lambda(b,1/b)\rd \mu(b)$ can be calculated analytically.

Asymptotic results for empirical tail copulas are derived under the following assumptions:
\begin{enumerate}[label=\arabic*)]
    \item\label{item:assumption:non:zero:Lambda} $\Lambda(\cdot)=\Lambda(\cdot;C)\not \equiv 0$;
    \item\label{item:assumption:bounds:convergence:empirical:copula} there exists a function $A:\IR_{+}\rightarrow \IR_{+}$ such that
    $\lim_{t\rightarrow \infty}A(t)=0$ and that $|\Lambda(u,v)-tC(u/t,v/t)|=\operatorname{O}(A(t))$, $t\to\infty$,
    locally uniformly for $(u,v) \in \IR_{+}^2$;
    \item\label{item:assumption:sequence:k} there exists a sequence $k=k(n) \in \IN$ such that $k\rightarrow \infty$, $k=\operatorname{o}(n)$ and $\sqrt{k}A\left(n/k\right)\rightarrow 0$; and
    \item\label{item:assumption:derivative} the partial derivative $\mathrm{D}_j\Lambda$ exists and is continuous on $\{(x_1,x_2)\in\IR_{+}^2:0<x_j<\infty\}$ for $j=1,2$;
\end{enumerate}
see~\citet{bucher2013multiplier} and~\citet{bormann2020detecting} for more details.
Under these assumptions, asymptotic results for the estimators in~\eqref{eq:estimators:mtcm},~\eqref{eq:estimators:mu:atcm:known:margin} and~\eqref{eq:estimators:mu:atcm:unknown:margin} can be derived as follows.

\begin{theorem}[Asymptotic normality for the MTCM]\label{theorem:asymptotic:normality:mtcm}
Suppose that the map $b \mapsto \Lambda(b,1/b)$ admits a unique maximum attained at $b^\ast \in [1/L,L]$ for some $L>1$.
If Assumptions~\ref{item:assumption:non:zero:Lambda},~\ref{item:assumption:bounds:convergence:empirical:copula} and~\ref{item:assumption:sequence:k} hold, then
\begin{align*}
\sqrt{k}\left({\tilde \lambda}^{\ast[n,k,L]}-\lambda^\ast(\Lambda)\right)\darrow
\N(0,\tilde \tau_{\Lambda}(b^\ast)),\quad \tilde \tau_{\Lambda}(b^\ast)=\Lambda\left(b^\ast,\frac{1}{b^\ast}\right).
\end{align*}
If, in addition, Assumption~\ref{item:assumption:derivative} is satisfied, then
\begin{align*}
\sqrt{k}\left({\hat \lambda}^{\ast[n,k,L]}-\lambda^\ast(\Lambda)\right)\darrow
\N(0,\hat \tau_{\Lambda}(b^\ast)),\quad \hat \tau_{\Lambda}(b^\ast)=\Var\left(\hat{\mathbb G}_\Lambda\left(b^\ast,\frac{1}{b^\ast}\right)\right),
\end{align*}
where the process $\hat {\mathbb G}_\Lambda$ is defined by
\begin{align}\label{eq:def:hat:G:process}
\hat{\mathbb G}_\Lambda(u,v)=\tilde{\mathbb G}_\Lambda(u,v)-\mathrm{D}_1 \Lambda(u,v)\tilde{\mathbb G}_\Lambda(u,\infty)-\mathrm{D}_2 \Lambda(u,v)\tilde{\mathbb G}_\Lambda(\infty,v),
\end{align}
with $\mathrm{D}_j\Lambda(u,v)$, $j=1,2$, being defined as $0$ if $u=0$ or $v=0$, and with $\tilde {\mathbb G}_\Lambda$ being a centered tight continuous Gaussian random field with the covariance structure given by $\E[\tilde{\mathbb G}_\Lambda(u,v)\tilde{\mathbb G}_\Lambda(u',v')]=\Lambda(u\wedge u',v\wedge v')$.
\end{theorem}

\begin{theorem}[Asymptotic normality for $\mu$-ATCMs]\label{theorem:asymptotic:normality:mu:atcm}
Suppose that the support of $\mu \in \mathcal M$ is contained in $[1/L,L]$ for some $L>1$.
If Assumptions~\ref{item:assumption:non:zero:Lambda},~\ref{item:assumption:bounds:convergence:empirical:copula} and~\ref{item:assumption:sequence:k} hold, then
\begin{align*}
\sqrt{k}\left(\tilde \lambda_{\mu}^{[n,k]}-\lambda_\mu(\Lambda)\right)\darrow
\N(0,\tilde \sigma_{\mu}^2(\Lambda)),
\end{align*}
where
\begin{align*}
\tilde \sigma_{\mu}^2(\Lambda)=\frac{\int_{(0,\infty)}\!\int_{(0,\infty)} \Lambda(\min(b,b'),\min(1/b,1/b'))\rd \mu(b)\!\rd \mu(b')}{\left(\int_{(0,\infty)}\overline \Lambda(b,1/b)\rd \mu(b)\right)^2}.
\end{align*}
If, in addition, Assumption~\ref{item:assumption:derivative} is satisfied, then\begin{align*}
\sqrt{k}\left(\hat \lambda_{\mu}^{[n,k]}-\lambda_\mu(\Lambda)\right)\darrow
\N(0,\hat \sigma_{\mu}^2(\Lambda)),
\end{align*}
where
\begin{align*}
\hat \sigma_{\mu}^2(\Lambda)=\frac{\int_{(0,\infty)}\!\int_{(0,\infty)} \E[\hat {\mathbb G}_{\Lambda}(b,1/b)\hat {\mathbb G}_{\Lambda}(b',1/b')]\rd \mu(b)\!\rd \mu(b')}{\left(\int_{(0,\infty)}\overline \Lambda(b,1/b)\rd \mu(b)\right)^2},
\end{align*}
where $\hat {\mathbb G}$ is as defined in~\eqref{eq:def:hat:G:process}.
\end{theorem}

\begin{remark}
In Theorems~\ref{theorem:asymptotic:normality:mtcm} and~\ref{theorem:asymptotic:normality:mu:atcm}, the estimators~\eqref{eq:estimators:mtcm},~\eqref{eq:estimators:mu:atcm:known:margin} and~\eqref{eq:estimators:mu:atcm:unknown:margin} are constructed based on the restricted part of the tail copula $\left\{\Lambda\left(b,1/b\right): b \in \left[1/L,L\right]\right\}$.
This restriction simplifies the technical discussion in the proof of the above theorems, and more advanced analysis is left for future research.
\end{remark}

\subsection*{Proofs of Theorems~\ref{theorem:asymptotic:normality:mtcm} and~\ref{theorem:asymptotic:normality:mu:atcm}}

Let $\mathcal B_\infty(\bar\IR_{+}^2)$ denote the space of all $\IR$-valued functions on $\bar\IR_{+}^2$ which are locally uniformly bounded on every compact subset of $[0,\infty]^2\backslash\{(0,0)\}$.
The set $\mathcal B_\infty(\bar\IR_{+}^2)$ is a complete metric space with the metric
\begin{align*}
d(f,g)=\sum_{l=1}^\infty \frac{1\wedge ||f-g||_{T_l}}{2^l},\quad f,\,g \in \mathcal B_\infty(\bar\IR_{+}^2),
\end{align*}
where $T_l=[0,l]^2\cup[0,l]\times \{\infty\}\cup \{\infty\}\cup[0,l]$ and $||f||_{T_l}=\sup_{x \in T_l}|f(x)|$; see \citet[Chapter~1.6]{van1996weak}.
Note that a sequence in $\mathcal B_\infty(\bar\IR_{+}^2)$ converges with respect to this metric if and only if it converges uniformly on every $T_l$.
Let $\rightsquigarrow$ denote weak convergence in the sense of Hoffmann-J{\o}rgensen; see the aforementioned reference for details.
Then the following asymptotic results hold for empirical tail copulas; see \citet[Lemma~2.1 and Theorem~2.2]{bucher2013multiplier}.

\begin{proposition}
\label{prop:asymptotic:normality:tail copula:known:margins}
Suppose that Assumptions~\ref{item:assumption:non:zero:Lambda},~\ref{item:assumption:bounds:convergence:empirical:copula} and~\ref{item:assumption:sequence:k} hold.
Then, as $n\rightarrow \infty$, we have that
\begin{align}\label{eq:result:weak:convergence:empirical:tail:copula:known:margins}
\sqrt{k}\left(\tilde \Lambda^{[n,k]}(u,v)-\Lambda(u,v)\right)\rightsquigarrow\tilde {\mathbb G}_{\Lambda}(u,v),
\end{align}
in $\mathcal B_\infty(\bar\IR_{+}^2)$, where $\tilde {\mathbb G}_\Lambda$ is a centered tight continuous Gaussian random field with covariance structure
$\E[\tilde{\mathbb G}_\Lambda(u,v)\tilde{\mathbb G}_\Lambda(u',v')]=\Lambda(u\wedge u',v\wedge v')$.
\end{proposition}

\begin{proposition}
\label{prop:asymptotic:normality:tail copula:unknown:margins}
Suppose that Assumptions~\ref{item:assumption:non:zero:Lambda},~\ref{item:assumption:bounds:convergence:empirical:copula},~\ref{item:assumption:sequence:k} and~\ref{item:assumption:derivative} hold.
Then, as $n\rightarrow \infty$, we have that
\begin{align}\label{eq:result:weak:convergence:empirical:tail:copula:unknown:margins}
\sqrt{k}\left(\hat \Lambda^{[n,k]}(u,v)-\Lambda(u,v)\right)\rightsquigarrow\hat {\mathbb G}_{\Lambda}(u,v),
\end{align}
in $\mathcal B_\infty(\bar\IR_{+}^2)$, where $\hat {\mathbb G}_\Lambda$ is a centered tight continuous Gaussian random field represented by
\begin{align*}
\hat{\mathbb G}_\Lambda(u,v)=\tilde{\mathbb G}_\Lambda(u,v)-\mathrm{D}_1 \Lambda(u,v)\tilde{\mathbb G}_\Lambda(u,\infty)-\mathrm{D}_2 \Lambda(u,v)\tilde{\mathbb G}_\Lambda(\infty,v),
\end{align*}
and $\mathrm{D}_j\Lambda(u,v)$, $j=1,2$, is defined as $0$ if $u=0$ or $v=0$.
\end{proposition}

By~\citet[Theorem~1.6.1]{van1996weak}, the above weak convergence results hold in the set of uniformly bounded functions $l^\infty(T_L)$ on $T_L=[0,L]^2\cup[0,L]\times \{\infty\}\cup \{\infty\}\times [0,L]$.
Therefore,~\eqref{eq:result:weak:convergence:empirical:tail:copula:known:margins} and~\eqref{eq:result:weak:convergence:empirical:tail:copula:unknown:margins} hold in $l^\infty([1/L,L])$ with $(u,v) \in [1/L,L]^2$ replaced by $(b,1/b)$, $b \in [1/L,L]$, and the domain of all the appearing functions is regarded as $[1/L,L]$; this statement is a consequence of the continuous mapping theorem~\citep[Theorem~1.3.6]{van1996weak} by taking $g:l^\infty(T_L)\rightarrow l^\infty([1/L,L])$ as $g(f)=f(b,1/b)$, $b \in [1/L,L]$.

Now Theorem~\ref{theorem:asymptotic:normality:mtcm} is derived directly from Propositions~\ref{prop:asymptotic:normality:tail copula:known:margins} and~\ref{prop:asymptotic:normality:tail copula:unknown:margins} and the extended Delta Method of \citet[Theorem~2.2]{carcamo2020directional} under the assumption that the unique maximizer $b^\ast$ of $b \rightarrow \Lambda(b,1/b)$ is in $[1/L,L]$.
Namely, for ${\tilde \lambda}^{\ast[n,k,L]}$, we have that
\begin{align*}
\sqrt{k}\left({\tilde \lambda}^{\ast[n,k,L]}-\lambda^\ast(\Lambda)\right)
&=\sqrt{k}\left( \max_{b \in [1/L,L]}{\tilde \Lambda}^{[n,k]}\left(b,\frac{1}{b}\right)- \max_{b \in [1/L,L]}{\Lambda}\left(b,\frac{1}{b}\right)\right)\\
&\darrow
\sup
\left\{
\tilde{\mathbb G}_\Lambda\left(b,\frac{1}{b}\right)
: b \in \left[\frac{1}{L},L\right]
\text{ such that }
\Lambda\left(b,\frac{1}{b}\right)=
\Lambda\left(b^\ast,\frac{1}{b^\ast}\right)
\right\}
=\tilde{\mathbb G}_\Lambda\left(b^\ast,\frac{1}{b^\ast}\right),
\end{align*}
where the limit is obtained by~\citet[Corollary 2.3]{carcamo2020directional}.
The case of ${\hat \lambda}^{\ast[n,k,L]}$ is shown analogously.

Next, the weak convergence results in Theorem~\ref{theorem:asymptotic:normality:mu:atcm} hold directly by Propositions~\ref{prop:asymptotic:normality:tail copula:known:margins} and~\ref{prop:asymptotic:normality:tail copula:unknown:margins} and the continuous mapping theorem.
Namely, since the map $\phi^{[L]}: l^\infty([1/L,L])\rightarrow \IR$ defined by
\begin{align*}
\phi^{[L]}(f)
&=\frac{\int_{(0,\infty)}f(b)\rd \mu(b)}{\int_{(0,\infty)}\overline \Lambda(b,1/b)\rd \mu(b)}
\end{align*}
is linear, we have, for ${\tilde \lambda}_\mu^{[n,k]}$, that
\begin{align*}
\sqrt{k}\left({\tilde \lambda}_{\mu}^{[n,k]}-\lambda_\mu(\Lambda)
\right)
&=\sqrt{k}\left({\tilde \lambda}_{\mu}^{[n,k]}-\lambda_\mu(\Lambda)\right)\\
&=\frac{1}{\int_{(0,\infty)}\overline \Lambda(b,1/b)\rd \mu(b)}
\int_{(0,\infty)}
\sqrt{k}\left(
{\tilde \Lambda}^{[n,k]}\left(b,\frac{1}{b}\right)- \Lambda\left(b,\frac{1}{b}\right)
\right)
\rd \mu(b)\\
& \darrow
\frac{1}{\int_{(0,\infty)}\overline \Lambda(b,1/b)\rd \mu(b)}
\int_{(0,\infty)}
\tilde{\mathbb G}_\Lambda\left(b,\frac{1}{b}\right)
\rd \mu(b).
\end{align*}
By~\citet[Lemma~3.9.8]{van1996weak}, the limit is normally distributed with mean
\begin{align*}
\frac{
\int_{(0,\infty)}
\E\left[\tilde{\mathbb G}_\Lambda\left(b,\frac{1}{b}\right)\right]
\rd \mu(b)}
{
\int_{(0,\infty)}\overline \Lambda(b,1/b)\rd \mu(b)
}
=0
\end{align*}
and variance
\begin{align*}
\frac{
\int_{(0,\infty)}\int_{(0,\infty)}
\E\left[\tilde{\mathbb G}_\Lambda\left(b,\frac{1}{b}\right)
\tilde{\mathbb G}_\Lambda\left(b',\frac{1}{b'}\right)
\right]
\rd \mu(b)\rd \mu(b')
}
{
	\left(\int_{(0,\infty)}\overline \Lambda(b,1/b)\rd \mu(b)\right)^2
}
=\frac{
\int_{(0,\infty)}\int_{(0,\infty)}
\Lambda\left(\min(b,b'),\min\left(\frac{1}{b},\frac{1}{b'}\right)\right)
\rd \mu(b)\rd \mu(b')
}
{
	\left(\int_{(0,\infty)}\overline \Lambda(b,1/b)\rd \mu(b)\right)^2
}.
\end{align*}
The case of ${\hat \lambda}_\mu^{[n,k]}$ is shown analogously.
\end{document}